\documentclass[a4paper,12pt]{article}
\usepackage{dsfont}
\usepackage{amssymb}
\usepackage{latexsym}
\usepackage{amsmath}
\usepackage{color}
\usepackage{comment}
\usepackage{amsthm}% for \proof,\newtheorem
\usepackage[dvips]{graphicx}% get graphs
\usepackage{layout} % for layout 
\usepackage{ulem}
\usepackage{enumerate}%\begin{enumerate}[{(}i{)}]
\usepackage{bm}
\usepackage{bbm} %\mathbbm{ABCabc\nabla}
\usepackage[bbgreekl]{mathbbol} %\mathbb{abcABC\nabla\bbalpha\bbbeta\Delta  }
\usepackage{authblk}
\newif\ifrs
\rstrue
\ifrs \usepackage{mathrsfs} \fi  % Use \mathscr{*}
\newif\ifcol
\colfalse %%% Deactivate "notice"

\newtheorem{theorem*}{Theorem}[section]
\newtheorem{note*}[theorem*]{Note}
\newtheorem{lemma*}[theorem*]{Lemma}
\newtheorem{definition*}[theorem*]{Definition}
\newtheorem{proposition*}[theorem*]{Proposition}
\newtheorem{corollary*}[theorem*]{Corollary}
\newtheorem{remark*}[theorem*]{Remark}
\newtheorem{example*}[theorem*]{Example}
\numberwithin{equation}{section}
\newcommand{\colred}{\color{black}}%{\color[rgb]{0.8,0,0}}
\newif\ifcol
\colfalse
\ifcol
\newcommand{\colorr}{\color[rgb]{0.8,0,0}}
\newcommand{\colorg}{\color[rgb]{0,0.5,0}}
\newcommand{\colorb}{\color[rgb]{0,0,0.8}}

\newcommand{\colorn}{\color[rgb]{1,1,1}}
\newcommand{\coloroy}{\color[rgb]{1,0.95,0}}

\else
\newcommand{\colorb}{\color{black}}% {\color[rgb]{0,0,0.8}}
\newcommand{\colorr}{\color{black}}% {{\color[rgb]{0.8,0,0}}
\newcommand{\colorg}{\color{black}}% {\color[rgb]{0,0.5,0}}
\newcommand{\colorn}{\color{black}}% {\color[rgb]{1,1,1}}
\newcommand{\coloroy}{\color{black}}% {\color[rgb]{1,0.95,0}}
\fi
\excludecomment{en-text}
\includecomment{jp-text}
\includecomment{comment}
\newtheorem{theorem}{Theorem}[section]

\newtheorem{lemma}[theorem*]{Lemma}

\newtheorem{proposition}[theorem*]{Proposition}

\newtheorem{remark}[theorem*]{Remark}

\numberwithin{equation}{section}
\def\koko{{\coloroy{koko}}}
\def\bd{\begin{description}}
\def\ed{\end{description}}

\def\D2{\bbD_{2,\infty-}}

\def\C{{\bf C}}
\def\D{{\bf D}}

\def\calb{{\cal B}}
\def\calc{{\cal C}}

\def\calf{{\cal F}}
\def\calg{{\cal G}}
\def\calh{{\cal H}}

\def\calp{{\cal P}}

\def\calt{{\cal T}}

\def\sfr{{\sf r}}
\def\sfk{{\sf k}}

\def\yleq{\>\leq\>}

\def\ds{\displaystyle}

\def\yeq{\>=\>}
\def\sfm{{\sf m}}

\def\sfd{{\sf d}}
\def\sfp{{\sf p}}

\def\simleq{\ \raisebox{-.7ex}{$\stackrel{{\textstyle <}}{\sim}$}\ }

\def\ep{\epsilon}
\def\half{\frac{1}{2}}

\def\inclusion{\hookrightarrow}
\def\cadlag{c\`adl\`ag\ }
\def\up{\uparrow}
\def\down{\downarrow}

\def\y{\vspace*{3mm}\\}
\def\halflineskip{\vspace*{3mm}}
\def\nn{\nonumber}
\def\be{\begin{equation}}
\def\ee{\end{equation}}
\def\bea{\begin{eqnarray}}
\def\eea{\end{eqnarray}}
\def\beas{\begin{eqnarray*}}
\def\eeas{\end{eqnarray*}}
\def\bi{\begin{itemize}}
\def\ei{\end{itemize}}
\def\im{\item}
\def\bd{\begin{description}}
\def\ed{\end{description}}
\def\l{\left}
\def\r{\right}
\def\dotc{\stackrel{\circ}{C}}

\newcommand{\bbB}{{\mathbb B}}

\newcommand{\bbD}{{\mathbb D}}

\newcommand{\bbF}{{\mathbb F}}
\newcommand{\bbG}{{\mathbb G}}
\newcommand{\bbH}{{\mathbb H}}

\newcommand{\bbN}{{\mathbb N}}

\newcommand{\bbQ}{{\mathbb Q}}
\newcommand{\bbR}{{\mathbb R}}

\newcommand{\bbT}{{\mathbb T}}
\newcommand{\bbU}{{\mathbb U}}
\newcommand{\bbV}{{\mathbb V}}

\newcommand{\bbY}{{\mathbb Y}}
\newcommand{\bbZ}{{\mathbb Z}}

\begin{document}
%\pagecolor{black}\color{white} %white on black

%%%%%%%%%%%%%%%%%%%%%%%%%%%%%%%%%%%%%%%%%
%%%%%%%%%%%%%%%%%%%%%%%%%%%%%%%%%%%%%%%%%
\title{%Expansion of the Asymptotically Conditionally Normal Law %arXiv
%Expansion of the Asymptotically Conditionally Normal Law %ISM Research Memo
%Partial quasi likelihood analysis and partial mixing 
Partial quasi likelihood analysis%Submit
\footnote{
This work was in part supported by 
CREST JPMJCR14D7 Japan Science and Technology Agency; 
Japan Society for the Promotion of Science Grants-in-Aid for Scientific Research 
No. 17H01702 (Scientific Research) 
%No. 24340015 (Scientific Research), 
%Nos. 24650148 and 
%No. 26540011 (Challenging Exploratory Research); \koko
%the Global COE program ``The Research and Training Center for New Development in Mathematics'' of the Graduate School of Mathematical Sciences, University of Tokyo; 
%NS Solutions Corporation; 
and by a Cooperative Research Program of the Institute of Statistical Mathematics. 
%
%The main parts in this paper were presented at 
%International conference ``Statistique Asymptotique des Processus Stochastiques VII'', Universit\'e du Maine, Le Mans, March 16-19, 2009, 
%MSJ Spring Meeting 2010, March 24-27, 2010, Keio University, Mathematical Society of Japan, 
%and 
%International conference ``DYNSTOCH Meeting 2010'', Angers, June 16-19, 2010. 
%The author thanks to the organizers of the meetings for opportunities of the talks. 
}
}
%%%% Author %%%%%%%%%%%%%%%%%%%%%%%
%%%%%%%%%%%%%%%%%%%%%%%%%%%%%%%
\author{Nakahiro Yoshida}
\affil{Graduate School of Mathematical Sciences, University of Tokyo
\footnote{Graduate School of Mathematical Sciences, University of Tokyo: 3-8-1 Komaba, Meguro-ku, Tokyo 153-8914, Japan. e-mail: nakahiro@ms.u-tokyo.ac.jp}}
\affil{CREST, Japan Science and Technology Agency%\footnote{}
}
%%%% Authors %%%%%%%%%%%%%%%%%%%%%%%
%%%%%%%%%%%%%%%%%%%%%%%%%%%%%%%
%\author[1,3]{Masayuki Uchida}
%\author[2,3]{Nakahiro Yoshida}
%\affil[1]{Graduate School of Engineering Science, Osaka University
%\footnote{Graduate School of Engineering Science, Osaka University: Toyonaka, Osaka 560-8531, Japan}
%        }
%\affil[2]{Graduate School of Mathematical Sciences, University of Tokyo
%\footnote{Graduate School of Mathematical Sciences, University of Tokyo: 3-8-1 Komaba, Meguro-ku, Tokyo 153-8914, Japan. e-mail: nakahiro@ms.u-tokyo.ac.jp}
%        }
%\affil[3]{CREST, Japan Science and Technology Agency
%%\footnote{}
%        }
%%%% Date %%%%%%%%%%%%%%%%%%%%%%%%
%%%%%%%%%%%%%%%%%%%%%%%%%%%%%%%
\date{December 29, 2017\\
%%Revised February 25, 2012
}
%%%%%%%%%%%%%%%%%%%%%%%%%%%%%%%
\maketitle
%%%%%%%%%%%%%%%%%%%%%%%%%%%%%%%
%%%%%%%%%%%%%%%%%%%%%%%%%%%%%%%
\ \\
{\it Summary} \ 
The quasi likelihood analysis is generalized to the partial quasi likelihood analysis. 
Limit theorems for the quasi likelihood estimators, especially the quasi Bayesian estimator, are derived 
in the situation where existence of a slow mixing component prohibits the Rosenthal type inequality from applying to the derivation of 
%applying a that would derive a 
the polynomial type large deviation inequality for the statistical random field. 
%without assuming existence of the moment of the statistical random field. 
%in the situation where the statistical random field does not admit high order of moments but admits them under conditioning. 
We give two illustrative examples. %estimation of a diffusion process having a component with a slow mixing rate.  

%%%%%%%%%%%%%%%%%%%%%%%%%%%%%%%%%%%%%%%%%%%%%%%%%%%%%%%%%%%%
%%%%%%%%%%%%%%%%%%%%%%%%%%%%%%%%%%%%%%%%%%%%%%%%%%%%%%%%%%%%
%%%%%%%%%%%%%%%%%%%%%%%%%%%%%%%%%%%%%%%%%%%%%%%%%%%%%%%%%%%%
%%%%%%%%%%%%%%%%%%%%%%%%%%%%%%%%%%%%%%%%%%%%%%%%%%%%%%%%%%%%
%%%%%%%%%%%%%%%%%%%%%%%%%%%%%%%%%%%%%%%%%%%%%%%%%%%%%%%%%%%%
\section{Introduction}

The Ibragimov-Has'minskii theory enhanced 
the asymptotic decision theory by Le Cam and H\'ajek by 
convergence of the likelihood ratio random field,  
and was programed by Kutoyants to statistical inference for semimartingales. 
The core of the theory is the large deviation inequality for the 
associated likelihood ratio random field.  
Asymptotic properties of the likelihood estimators are deduced from those of the likelihood ratio random field. 
Precise estimates of the tail probability and hence convergence of moments of the estimators follow 
in a unified manner once such a strong mode of convergence of the likelihood ratio random field is established. 
For details, see Ibragimov and Has'minskii \cite{IbragimovHascprimeminskiui1972,IbragimovHascprimeminskiui1973,IbragimovHascprimeminskiui1981} and 
Kutoyants \cite{Kutoyants1984,Kutoyants1994,Kutoyants1998,Kutoyants2004}. 

The quasi likelihood analysis
 (QLA) descended from the Ibragimov-Has'minskii-Kutoyants program.\noindent
 \footnote{The QLA is not in the sense of Robert Wedderburn. Since exact likelihood function can rarely be assumed in inference for 
discretely sampled continuous time processes, quasi likelihood functions are quite often used there. 
Further, the word ``QLA'' also implies a new framework of inferential theory for stochastic processes within which 
the polynomial type large deviation is easily available today and plays an essential role in the theory. } 
In Yoshida \cite{yoshida2011polynomial}, 
it was showed that a polynomial type large deviation (PLD) inequality universally follows from 
certain separation of the random field, such as the local asymptotic quadraticity of the random field, 
and $L_p$ estimates of easily tractable random variables. 
Since the PLD inequality is no longer a bottleneck of the program, the QLA applies to various complex random fields. % associated with nonlinear stochastic models and 
\begin{en-text}
It showed that certain separation of the random field such as the local asymptotic quadraticity of the random field 
is sufficient for the polynomial type large deviation inequality and that 
the large deviation estimate of the random field %on the tangent space of the parameter space 
can be reduced to $L_p$ estimates of easily tractable random variables. 
This approach has merit 
Various statistical models enjoyed merits from this approach.  
\end{en-text}

The QLA is a framework of statistical inference for stochastic processes. 
It features the polynomial type large deviation of the quasi likelihood random field. 
%Various inferential problems enjoyed merit from this scheme. 
Through QLA, one can systematically derive 
limit theorems and precise tail probability estimates of the associated QLA estimators such as quasi maximum likelihood estimator (QMLE), 
quasi Bayesian estimator (QBE) and various adaptive estimators. 
The importance of such precise estimates of tail probability is well recognized in asymptotic decision theory, 
prediction, theory of information criteria for model selection, asymptotic expansion, etc. 
The QLA is rapidly expanding the range of its applications: for example, 
sampled ergodic diffusion processes (Yoshida \cite{yoshida2011polynomial}), 
%misspecified diffusion processes (Uchida and Yoshida \cite{UchidaYoshida2011}, 
contrast-based information criterion for diffusion processes (Uchida \cite{Uchida2010}), 
approximate self-weighted LAD estimation of discretely observed ergodic Ornstein-Uhlenbeck processes (Masuda \cite{Masuda2010a}), 
jump diffusion processes Ogihara and Yoshida(\cite{OgiharaYoshida2011}), 
adaptive estimation for diffusion processes (Uchida and Yoshida \cite{UchidaYoshida2012Adaptive}), 
adaptive Bayes type estimators for ergodic diffusion processes (Uchida and Yoshida \cite{uchida2014adaptive}), 
%estimation of volatility of a sampled semimartingale in finite time horizon is typical non-ergodic statistics. 
asymptotic properties of the QLA estimators for volatility in regular sampling of finite time horizon (Uchida and Yoshida \cite{UchidaYoshida2013})  
and in non-synchronous sampling (Ogihara and Yoshida \cite{ogihara2014quasi}), 
Gaussian quasi-likelihood random fields for ergodic L\'evy driven SDE (Masuda \cite{masuda2013convergence}), 
hybrid multi-step estimators (Kamatani and Uchida  \cite{KamataniUchida2014}), 
parametric estimation of L{\'e}vy processes (Masuda \cite{masuda2015parametric}),
ergodic point processes for limit order book (Clinet and Yoshida \cite{clinet2015statistical}), 
a non-ergodic point process regression model (Ogihara and Yoshida \cite{ogihara2015quasi}), 
threshold estimation for stochastic processes with small noise (Shimizu \cite{shimizu2015threshold}), 
AIC for non-concave penalized likelihood method (Umezu et al. \cite{umezu2015aic}), 
%a spot volatility information criterion (sVIC) for volatility model selection (Uchida and Yoshida \cite{uchida2016model}), 
Schwarz type model comparison for LAQ models (Eguchi and Masuda \cite{eguchi2016schwarz}), 
adaptive Bayes estimators and hybrid estimators for small diffusion processes based on sampled data (Nomura and Uchida \cite{nomura2016adaptive}), 
moment convergence of regularized least-squares estimator for linear regression model (Shimizu \cite{shimizu2017moment}),  
moment convergence in regularized estimation under multiple and mixed-rates asymptotics (Masuda and Shimizu \cite{masuda2017moment}), 
asymptotic expansion in quasi likelihood analysis for volatility (Yoshida \cite{Yoshida2017asymptoticexpansion}) among others.

As already mentioned, the PLD inequality is the key to the QLA. 
Once a PLD inequality is established, we can obtain a very strong mode of convergence of the random field %$\bbZ_T$. 
and the associated estimators. 
%as a result, 
%that gives convergence of moments of the QLA estimators as well as limit theorems for the QBE. 
However, in the present theory, boundedness of high order of moments of functionals 
is assumed. %associated with $\bbH_T$ is necessary 
% to apply the QLA theory. 
On the other hand, for example, if the statistical model has a component with a slow mixing rate, 
the Rosenthal inequality does not serve to validate the boundedness of moments of very high order. 
How do QMLE and QBE behave in such a situation? 
This question motivates us to introduce the partial quasi likelihood analysis (PQLA). 

The aim of this short note is to formulate the PQLA and to exemplify it. 
The basic idea is conditioning by partial information. %a sub $\sigma$-field $\calc$ of $\calf$. 
Easy to understand is a situation where there are two components $(L,U)$ of stochastic processes and 
$U$ has a fast mixing rate but $L$ has a slow mixing rate. 
Suppose that the Rosenthal inequality may control the moments of a functional of $U$ but 
cannot control the moments of a functional of $L$. 
In this situation, we cannot apply the present QLA theory or the way of derivation of the PLD inequality 
to the random fields expressed by $U$ and $L$. 
However, if there is a partial mixing structure in that $U$ possesses a very good mixing rate conditionally on $L$, 
then we can apply a conditional version of the QLA theory for given $L$. 
Even if $L$ has a bad mixing rate and its temporal impact on the system is unbounded, 
there is a possibility that we can recover limit theorems for the QLA estimators. 
Technically, a method of truncation is essential to detach the slow mixing component's effects %which potentially have unbounded impacts, 
from the main body of the randomness. 

Partial QLA naturally emerges in the structure of the partial mixing. 
The notion of partial mixing was used in Yoshida \cite{yoshida2004partial} to derive asymptotic expansion of 
the distribution of an additive functional of the conditional $\ep$-Markov process admitting a component with long-range dependency. 
%In other words, without truncation, we cannot apply Rosenthal type of inequalities under the conditional probability. 

The organization of this note is as follows. 
Section \ref{170906-10} presents a frame of the partial quasi likelihood analysis. 
The asymptotic properties of the QMLE and QBE are provided there. 
The conditional polynomial type large deviation inequality is the key to the partial QLA. 
Section \ref{290924-1} gives a set of sufficient conditions for it. %the inequality. 
%Partial QLA naturally emerges in the structure of the partial mixing. 
A conditional version of a Rosenthal type inequality is stated in Section \ref{290924-2}. 
Section \ref{290924-5} illustrates a diffusion process having slow and fast mixing components. 
Statistics is ergodic in Section \ref{290924-5}, while a non-ergodic statistical problem will be discussed in Section \ref{290924-6}.

\section{Partial quasi likelihood analysis}\label{170906-10}

%\subsection{Localization}
%\subsection{Partial polynomial type large deviation inequality}

\subsection{Quasi likelihood analysis% on $(\Omega,\calf,P_\calc)$
}\label{290918-1}
Given a probability space $(\Omega,\calf,P)$, we consider a sequence of random fields 
$\bbH_T:\Omega\times\overline{\Theta}\to\bbR$, $T\in\bbT$, 
where $\bbT$ is a subset of $\bbR_+$ with $\sup\bbT=\infty$,  
$\Theta$ is a bounded domain in $\bbR^\sfp$ and $\overline{\Theta}$ is its closure. 
We assume that $\bbH_T$ is $\calf\otimes\bbB[\bbR^\sfp]$-measurable and that 
the mapping $\overline{\Theta}\ni\theta\mapsto\bbH_T(\omega,\theta)$ is continuous for every $\omega\in\Omega$. 
By convention, $\bbH_T(\omega,\theta)$ is simply denoted by $\bbH_T(\theta)$. 

The random field $\bbH_T$ serves like the log likelihood function 
%to carry out statistical inference for $\theta$ 
in the likelihood analysis, but does more.  
A measurable mapping $\hat{\theta}_T:\Omega\to\overline{\Theta}$ is called 
a quasi maximum likelihood estimator (QMLE) if 
\beas
\bbH_T(\hat{\theta}_T) &=& \max_{\theta\in\overline{\Theta}}\bbH_T(\theta)
\eeas
for all $\omega\in\Omega$. 
The mapping $\tilde{\theta}_T:\Omega\to\calc[\Theta]$, the convex hull of $\Theta$, is defined by 
\beas 
\tilde{\theta}_T &=& 
\bigg[\int_\Theta\exp\big(\bbH_T(\theta)\big)\varpi(\theta)d\theta\bigg]^{-1}\int_\Theta\theta\exp\big(\bbH_T(\theta)\big)\varpi(\theta)d\theta
\eeas
and called the quasi Bayesian estimator (QBE) with respect to the prior density $\varpi$. 
We assume $\varpi$ is continuous and satisfies $0<\inf_{\theta\in\Theta}\varpi(\theta)\leq\sup_{\theta\in\Theta}\varpi(\theta)<\infty$. 
We call these estimators together quasi likelihood estimators. 

%In analysis of the asymptotic properties of quasi likelihood estimators, 
The quasi likelihood analysis (QLA) is formulated with the random field 
\beas 
\bbZ_T(u) 
&=&
\exp\big(\bbH_T(\theta^*+a_Tu)-\bbH_T(\theta^*)\big)\quad(u\in\bbU_T)
\eeas
Here $\theta^*\in\Theta$ is the target value of $\theta$ in estimation and 
$\bbU_T=\{u\in\bbR^\sfp;\>\theta^\dagger_T(u)\in\Theta\}$, where $\theta^\dagger_T(u)=\theta^*+a_Tu$. 
The matrix $a_T\in\text{GL}(\bbR^\sfp)$ satisfies $a_T\to0$ as $T\to\infty$. 
It is possible to extend $\bbZ_T$ to $\bbR^\sfp$ so that the extension has a compact support and 
$\sup_{u\in\bbR^\sfp\setminus\bbU_T}\bbZ_T(u)\leq\max_{u\in\partial\bbU_T}\bbZ_T(u)$. 
We denote this extended random field by the same $\bbZ_T$. 
Let $\hat{C}=\{f\in C(\bbR^\sfp);\>\lim_{|u|\to\infty}f(u)=0\}$. 
Then $\bbZ_T\in\hat{C}$. 

Consider $\sigma$-fields $\calc$ and $\calg$ such that $\calc\subset\calg\subset\calf$. 
We introduce $\calc$-measurable variables $\Psi_T:\Omega\to\{0,1\}$. 
These functionals are helpful to localize QLA. 
\begin{en-text}
{\colorg Remove: 
We will assume that $(\Omega,\calf)$ is standard. Then a regular conditional distribution $P_\calc$ exists. 
However, it is possible to remove this assumption because every time $P_\calc$ appears, 
we can rewrite it by a conditional measure on a good function space in all situations we meet in this paper. }
\end{en-text}

\subsection{Quasi maximum likelihood estimator}
Let $L$ be a positive constant. 
We start with the so-called polynomial type large deviation inequality, which 
plays an essential role in the theory of QLA as in \cite{yoshida2011polynomial}. 
%In what follows, we will use some of them in typical settings but 
%more subdivided presentation is possible as in \cite{yoshida2011polynomial}, as a matter of fact. 
Let $\bbV_T(r)=\{u\in\bbU_T;\>|r|\geq r\}$. 
Let $B_{c,T}=\{u\in\bbR^\sfp;\> |u|<c,\>\theta^\dagger_T(u)\in\Theta\}$ for $c>0$. 
The modulus of continuity of $\log\bbZ_T$ is 
\beas 
w_T(\delta,c) &=& \sup\bigg\{\big|\log\bbZ_T(u_2)-\log\bbZ_T(u_1)\big|;\>
u_1,u_2\in B_{c,T},\>|u_2-u_1|\leq\delta\bigg\}.
\eeas

Let $W_T(\delta,c,\ep)=\big\{w_T(\delta,c)>\ep\big\}$ and let 
\beas 
S_T(r,\ep) &=& 
\bigg\{\sup_{u\in\bbV_T(r)}\bbZ_T(u)\geq \ep \bigg\}. 
\eeas
Let $T_0>0$. 
Let ${\mathfrak T}$ be the set of sequences $(T_n)_{n\in\bbN}$ of numbers in $\bbT$ 
such that $T_n\geq T_0$ for all $n\in\bbN$ and 
$\lim_{n\to\infty}T_n=\infty$. 
{\colred Let $(\Psi_T)_{T\in\bbT}$ be a sequence of $[0,1]$-valued $\calc$-measurable 
random variables. }

\bd
\im[[A1\!\!]] %Conditional PLD under $\Psi_T\cdot P_\calc$ 
There exists a sequence of positive $\calc$-measurable random variables $(\ep(r))_{r\in\bbN}$ such that $\lim_{r\to\infty}\ep(r)=0$ a.s. and that 
\beas 
\sup_{T\in\calt}\sup_{r\in\bbN}r^LP_\calc\big[S_T(r,\ep(r))\big]\Psi_T&<& \infty \quad a.s.
\eeas
for every $\calt\in{\mathfrak T}$. 
\ed
\bd
\im[[A1$^\flat$\!\!]] %Conditional PLD under $\Psi_T\cdot P_\calc$ 
For a sequence of positive numbers $(\ep(r))_{r\in\bbN}$ with $\lim_{r\to\infty}\ep(r)=0$ and 
a sequence of positive random variables $(\eta(r))_{r>0}$ with $\lim_{r\to\infty}\eta(r)=0$ a.s., it holds that 
\beas 
P_\calc\big[S_T(r,\ep(r))\big]\Psi_T&\leq& \eta(r) \quad a.s.\quad(T\in\calt,\>r\in\bbN)
\eeas
for every $\calt\in{\mathfrak T}$. 
\ed
\bd
\im[[A2\!\!]] 
$\ds \limsup_{T\to\infty,{\colorr T\in\calt}}P_\calc\big[W_T(\delta,c,\ep)\big]\Psi_T\to^P0$ as $\delta\downarrow0$ for every $\ep>0$, $c>0$ and $\calt\in{\mathfrak T}$. 
\ed

\begin{remark}\rm 
The estimate of modulus of continuity is used only countable times to prove tightness. 
\end{remark}

%Extension of $\bbZ_T$ to $\hat{C}$. 
We consider and its extension $(\overline{\Omega},\overline{\calf},\overline{P})$, that is, 
$\Omega\subset\overline{\Omega}$, $\calf\subset\overline{\calf}$ and $P=\overline{P}|_\calf$. 
Let $\bbZ(u)$ be a $\hat{C}$-valued random variable defined on 
$(\overline{\Omega},\overline{\calf},\overline{P})$. 
\footnote{Continuity is not necessary to assume as a matter of fact. Without it, the convergence of $\bbZ_T$ ensures the limit distribution is supported by 
$\hat{C}$; we may assume $\bbZ$ is a continuous process after modification. }

\bd
\im[[A3\!\!]] %Convergence of $\bbZ_T\Psi_T$ under $P_\calc$ on compacts
{\bf (i)} For any $k\in\bbN$, $u_i\in\bbR^\sfp$ $(i=1,...,k)$, $f\in C_b(\bbR^{k\sfp})$ and any bounded $\calg$-measurable random variable $Y$, 
\beas 
E_\calc\big[f\big((\bbZ_T(u_i))_{i=1,...,k}\big)\Psi_TY\big] &\to^P& \overline{E}_\calc\big[f\big((\bbZ(u_i))_{i=1,...,k}\big)Y\big]%\quad a.s.
\eeas
\bd
\im[(ii)] $\Psi_T\to^P1$. 
\ed
\ed

\bd
\im[[A4\!\!]] 
With probability one, there exists a unique element $\hat{u}\in\bbR^\sfp$ that maximizes $\bbZ$.  
\ed

Remark: {From} [A3] (i), we can remove $\Psi_T$ but keeping it explicitely is helpful in applications. 
We may assume $\hat{u}$ is $\overline{\calf}$-measurable; the given mapping $\hat{u}$ has a measurable version. 
The following theorems claim $\calc$-conditional $\calg$-stable convergence of $\bbZ_T$ and $\hat{u}_T=a_T^{-1}(\hat{\theta}_T-\theta^*)$.

\begin{theorem}\label{170902-1} 
%Let $L>\sfp$. 
Suppose that $[A1^\flat]$, $[A2]$ and $[A3]$ are satisfied. % and that $\Psi_T\to^P1$. 
Then 
\bea\label{170902-1.1}
E_\calc\big[F(\bbZ_T)Y\big] &\to^{P}& \overline{E}_\calc\big[F(\bbZ)Y\big]
\eea
for any $F\in C_b(\hat{C})$ and any bounded $\calg$-measurable random variable $Y$. 
In particular, 
\beas
\bbZ_T\to^{d_s(\calg)}\bbZ
\eeas
as $T\to\infty$. 
\end{theorem}

\begin{theorem}\label{170903-30}
%Let $L>\sfp$. 
Suppose that $[A1^\flat]$, $[A2]$, $[A3]$ and $[A4]$ are satisfied. % and that $\Psi_T\to^P1$. 
Then 
\bea\label{170902-2}
E_\calc\big[f(\hat{u}_T)Y\big]&\to^P&\overline{E}_\calc\big[f(\hat{u})Y\big]\quad 
\eea
for any $f\in C_b(\bbR^\sfp)$ and any bounded $\calg$-measurable random variable $Y$. 
In particular, 
\beas 
\hat{u}_T\to^{d_{s}(\calg)}\hat{u}
\eeas
as $T\to\infty$. 
\end{theorem}
\noindent
{\it Proof of Theorems \ref{170902-1} and \ref{170902-2}.} 
\def\bbz{\mathbb z} 
(a) 
\def\tbbZ{\widetilde{\mathbb Z}}
We may assume that $\|Y\|_\infty\leq1$. 
Let 
\beas 
\tbbZ_T(u)&=&
\left\{
\begin{array}{ll}
\bbZ_T(u)&\quad(\Psi_T=1)\y
e^{-|u|^2}&\quad(\Psi_T=0)
\end{array}
\right.
\eeas
In view of [A3] (ii), we may show 
\bea\label{170902-3}
E_\calc\big[F(\tbbZ_T)Y\big] &\to^{P}& \overline{E}_\calc\big[F(\bbZ)Y\big]
\eea
for $F\in C_b(\hat{C})$ in order to show (\ref{170902-1.1}). 
Then, by subsequence argument, it suffices to show that for any sequence $(T_n)$ 
with $0\leq T_1<T_2<\cdots\to\infty$, 
there exists a subsequence $(T_{n'})$ of $(T_n)$ such that (\ref{170902-3}) holds along $(T_{n'})$. 
For $k\in\bbN$, let ${\mathfrak G}_k$ be a countable subset of $C_b(\bbR^{k+1})$ that determines probability measures on $\bbR^{k+1}$. 
%Let ${\mathfrak G}=\cup_{k\in\bbN}{\mathfrak G}_k$. 

Let $P^{T_n}_\omega(d\bbz,dy)$ be a regular conditional distribution of $(\tbbZ_{T_n},Y)$ on $\hat{C}\times[-1,1]$ given $\calc$. 
Let $P_\omega$ be a regular conditional distribution of $(\bbZ,Y)$ given $\calc$. 
Moreover let 
\beas 
w(\delta,c,x)&=&\sup\big\{\big|x(u_2)-x(u_1)\big|;\>u_1,u_2\in B_{c,T},\>|u_2-u_1|\leq\delta\big\}
\eeas
for $x\in\hat{C}$, and 
let 
\beas 
\check{w}(c,x) &=& \sup\big\{\big|x(u)\big|;\>u\in\bbR^\sfp,\>|u|\geq c\big\}.
\eeas

According to $[A3]$ (ii) and $[A2]$, 
there exists a subsequence $(T_{n^{(1)}})$ of $(T_n)$ such that 
$\lim_{n^{(1)}\to\infty}\Psi_{T_{n^{(1)}}}=1$ {\colred a.s.} and that 
\beas 
\lim_{m\in\bbN,\>m\to\infty}
\limsup_{n^{(1)}\to\infty}P_\calc\big[W_{T_{n^{(1)}}}(m^{-1},j,k^{-1})\big]\Psi_{T_{n^{(1)}}}
&=&0\quad a.s.
\eeas
for all $j,k\in\bbN$. 
Moreover, from $[A1^\flat]$,  
for $k\in\bbN$, there exists an  $r_0>0$ such that 
$\ep(r)\leq k^{-1}$ for all $r\geq r_0$. 
Then 
\beas 
\lim_{r\to\infty}\limsup_{n^{(1)}\to\infty}P_\calc\bigg[S_{T_{n^{(1)}}}(r,k^{-1})\bigg]
&\leq&
\lim_{r\to\infty}\limsup_{n^{(1)}\to\infty}P_\calc\bigg[S_{T_{n^{(1)}}}(r,\ep(r))\bigg]\Psi_{T_{n^{(1)}}}
\>\leq\>\lim_{r\to\infty}\eta(r)\yeq0\quad a.s.
\eeas
Thus, 
thanks to [A1$^\flat$], [A2] and [A3], 
there exist an event $\Omega_0\in\calf$ with $P[\Omega_0]=1$ 
and a subsequence $(T_{\colred n'})$ of $(T_{\colred n^{(1)}})$ such that 
for any $\omega\in\Omega_0$, the following conditions hold: 
\bd
\im[(i)] 
$\ds
\lim_{m\to\infty}\limsup_{{\colred n'}\to\infty}P^{T_{\colred n'}}_\omega\bigg[\bigg\{(x,y);\ w(m^{-1},j,\log x)>k^{-1}\bigg\}\bigg]
\yeq0\quad(\forall j,k\in\bbN)
$
\im[(ii)] 
$\ds \lim_{j\to\infty}\limsup_{{\colred n'}\to\infty} P^{T_{{\colred n'}}}_\omega\bigg[\bigg\{(x,y);\ 
\check{w}(j,x)>k^{-1}\bigg\}\bigg] \yeq0\quad(\forall k\in\bbN)
$
\im[(iii)] For every $k\in\bbN$, 
\beas 
\int_{\hat{C}\times[-1,1]}g\big((x(u_i))_{i=1,...,k},y\big)P^{T_{{\colred n'}}}_\omega(dx,dy)
&\to& 
\int_{\hat{C}\times[-1,1]}g\big((x(u_i))_{i=1,...,k},y\big)P_\omega(dx,dy)%\quad(T\to\infty)
\eeas
as ${\colred n'}\to\infty$ for all $g\in{\mathfrak G}_k$. 
\im[(iv)] 
$\Psi_{T_{{\colred n'}}}\to1$ as ${\colred n'}\to\infty$. 
\im[(v)] 
$\ds 
\int_{\hat{C}\times[-1,1]}yP^{T_{{\colred n'}}}_\omega(dx,dy)
\yeq
\int_{\hat{C}\times[-1,1]}yP_\omega(dx,dy)%\quad(T\to\infty)
$ 
for all ${\colred n'}$. 
\ed
\halflineskip

For $\omega\in\Omega_0$, $\ep>0$ and $j,k\in\bbN$, there exist $m(k),j(k)\in\bbN$ such that 
\beas 
\sup_{{{\colred n'}}}P^{T_{{\colred n'}}}_\omega\big[A_{\omega,j,k}\big]
&<& 2^{-j-k-1}\ep
\eeas
and 
\beas 
\sup_{{{\colred n'}}} P^{T_{{\colred n'}}}_\omega\big[ B_{\omega,k}\big] ,
&<& 2^{-k-1}\ep
\eeas
where 
\beas 
A_{\omega,j,k} &=& 
\bigg \{(x,y);\ w(m(k)^{-1},j,\log\bbz)>k^{-1}\bigg\}
\eeas
and 
\beas 
B_{\omega,k}
&=&
\bigg\{(x,y);\ \check{w}(j(k),x)>k^{-1}\bigg\}, 
\eeas
respectively.
Let 
$A_\omega=\overline{\cap_{j,k\in\bbN}A_{\omega,j,k}\cap\cap_{k\in\bbN}B_{\omega,k}}$. 
Then $A_\omega$ is a compact set in $\hat{C}\times[-1,1]$ and 
$\sup_{{\colred n'}} P^{T_{{\colred n'}}}_\omega\big[A_\omega\big]\geq 1-\ep$. 
Therefore the family of probability measures 
$\{P^{T_{{\colred n'}}}_\omega\}_{{\colred n'}}$ is tight since $\tbbZ_{T_{{\colred n'}}}(0)=1$. 
{\colred Let $(n^\dagger)$ be any subsequence of $({\colred n'})$.}
Then 
there exist a subsequence ${\colred (n'')}$ of ${\colred (n^\dagger)}$, 
{\colred $(n'')$ depending on $\omega$,} and a probability measure $P^*_\omega$ on $\hat{C}\times[-1,1]$ 
such that 
$P^{T_{\colred n''}}_\omega\to P^*_\omega$ as ${\colred n''\to\infty}$. 
In particular, 
\beas 
\int_{\hat{C}\times[-1,1]}g\big((x(u_i))_{i=1,...,k},y\big)P^{T_{\colred n''}}_\omega(dx,dy)
&\to& 
\int_{\hat{C}\times[-1,1]}g\big((x(u_i))_{i=1,...,k},y\big)P^*_\omega(dx,dy)%\quad(T\to\infty)
\eeas
as ${\colred n''}\to\infty$ for every $\omega\in\Omega_0$ and every $g\in{\mathfrak G}_k$, $k\in\bbN$. 
Therefore 
\beas 
\int_{\hat{C}\times[-1,1]}g\big((x(u_i))_{i=1,...,k},y\big)P^*_\omega(dx,dy)
&=&
\int_{\hat{C}\times[-1,1]}g\big((x(u_i))_{i=1,...,k},y\big)P_\omega(dx,dy)
\eeas
for all $g\in{\mathfrak G}_k$, $k\in\bbN$. 
Since all finite-dimensional marginal distributions coincide, 
$P^*_\omega=P_\omega$.  
{\colred This implies $P^{T_{n^\dagger}}_\omega\to P^*_\omega$ as $n^\dagger\to\infty$, and hence}
\bea\label{170903-1}
P^{T_{n'}}_\omega\to P_\omega
\eea 
for every $\omega\in\Omega_0$. 
In particular, we obtain (\ref{170902-3}) along $(T_{n'})$, which gives  
Theorem \ref{170902-1}. \y
\noindent
(b) We may assume $0\leq Y\leq1$. 
Consider sequences $(T_n)$ and $(T_{n'})$ in Step (a). 
Let ${\mathfrak F}=\big\{\Delta_q(x);\>q\in\bbQ^\sfp\big\}$ with 
\beas 
\Delta_q(x) &=& 
(-1)\vee\bigg(\sup_{u\in R_q}x(u)-\sup_{u\in R_q^c}x(u)\bigg)\wedge 1\,
\eeas
where 
$R_q=\{u=(u_i)_{i=1,...,\sfp};\> u_i\leq q_i\>(i=1,...,\sfp)\}$ for $q=(q_i)_{i=1,...,\sfp}\in\bbQ^\sfp$. 
If $\omega\in\Omega_0$ satisfies 
$J_\omega:=\int_{\hat{C}\times[0,1]}yP_\omega(dx,dy)>0$, then 
the already obtained (\ref{170903-1}) yields 
\bea\label{170903-10}
\widetilde{P}_\omega^{T_{n'}} &\to& \widetilde{P}_\omega
\eea
as $n'\to\infty$, where 
\beas
\widetilde{P}_\omega^{T_{n'}}(B)=\int_{\hat{C}\times[0,1]}1_B(x)y\>dP_\omega^{T_{n'}}/J_\omega 
&\text{and} &
\widetilde{P}_\omega(B)=\int_{\hat{C}\times[0,1]}1_B(x)y\>dP_\omega/J_\omega
\eeas
for $B\in\bbB[\hat{C}]$. 
We notice that $\widetilde{P}_\omega^{T_{n'}}$ as well as $\widetilde{P}_\omega$ is a probability measure by (v) of Part (a). 
The convergence (\ref{170903-10}) gives 
\beas
\widetilde{P}_\omega\big[\Delta_q(x)>0\big]
&\leq&
\liminf_{n'\to\infty}\widetilde{P}^{T_{n'}}_\omega\big[\Delta_q(x)>0\big] 
\\&\leq&
\limsup_{n'\to\infty}\widetilde{P}^{T_{n'}}_\omega\big[\Delta_q(x)\geq0\big] 
\\&\leq& \widetilde{P}_\omega\big[\Delta_q(x)\geq0\big],
\eeas
for all $q\in\bbQ^\sfp$, 
and hence 
\bea\label{170903-11}
E_\calc\big[1_{\{\Delta_q(\bbZ)>0\}}Y\big]
&\leq&
\liminf_{n'\to\infty}E_\calc\big[1_{\{\Delta_q(\tbbZ_{T_{n'}})>0\}}Y\big]
\nn\\&\leq&
\limsup_{n'\to\infty}E_\calc\big[1_{\{\Delta_q(\tbbZ_{T_{n'}})\geq0\}}Y\big]
\nn\\&\leq&
\overline{E}_\calc\big[1_{\{\Delta_q(\bbZ)\geq0\}}Y\big]\quad a.s.
\eea
for all $q\in\bbQ^\sfp$% and all $\omega\in\Omega_0$
, since this is obvious when $J_\omega=0$. 
By definition, 
\beas 
\big\{\hat{u}_T\leq q\big\}\cap\{\Psi_T=1\} &\subset& \big\{\Delta_q(\tbbZ_T)\geq0\big\}
\eeas
and 
\beas 
\big\{\hat{u}_T\leq q\big\}\cup\{\Psi_T=0\} &\supset& \big\{\Delta_q(\tbbZ_T)>0\big\}
\eeas
Therefore (\ref{170903-11}) implies 
\bea\label{291008-1} 
\limsup_{n'\to\infty} E_\calc\big[1_{\{\hat{u}_{T_{n'}}\leq q\}}Y\big]
\nn&\leq&
\limsup_{n'\to\infty}E_\calc\big[1_{\{\Delta_q(\tbbZ_{T_{n'}})\geq0\}}Y\big]
\nn\\&\leq&
\overline{E}_\calc\big[1_{\{\Delta_q(\bbZ)\geq0\}}Y\big]
\nn\\&\leq&
\overline{E}_\calc\big[1_{\{\hat{u}\leq q\}}Y\big]\quad a.s.
\eea
for all $q\in\bbQ^\sfp$, and 
\bea\label{291008-2} 
\liminf_{n'\to\infty} E_\calc\big[1_{\{\hat{u}_{T_{n'}}\leq q\}}Y\big]
\nn&\geq&
\liminf_{n'\to\infty}E_\calc\big[1_{\{\Delta_q(\tbbZ_{T_{n'}})>0\}}Y\big]
\nn\\&\geq&
\overline{E}_\calc\big[1_{\{\Delta_q(\bbZ)>0\}}Y\big]
\nn\\&\geq&
\overline{E}_\calc\big[1_{\{\hat{u}<q\}}Y]\quad a.s.
\eea
for all $q\in\bbQ^\sfp$, where $(a_i)<(b_i)$ means $a_i<b_i$ for all $i=1,...,\sfp$, and we used uniqueness of $\hat{u}$ in the last part of each. 

Denote by $Q_\omega^{T_{n'}}$ [resp. $Q_\omega$] a regular conditional probability  of $(\hat{u}_{T_{n'}},Y)$ [resp. $(\hat{u},Y)$] 
given $\calc$. 
{From} (\ref{291008-1}) and (\ref{291008-2}), 
there exists $\Omega_1\in\calf$ with $P[\Omega_1]=1$ such that 
\beas 
I_\omega\>:=\>\int_{\bbR^\sfp\times[0,1]}y\>Q_\omega^{T_{n'}}(du,dy) &=& \int_{\bbR^\sfp\times[0,1]}y\>Q_\omega(du,dy)
\eeas
for all $n'$ and all $\omega\in\Omega_1$, and that 
\beas 
\int_{\bbR^\sfp\times[0,1]}1_{\{u< q\}}y\>Q_\omega(du,dy)
&\leq&
\liminf_{n'\to\infty}\int_{\bbR^\sfp\times[0,1]}1_{\{u\leq q\}}y\>Q_\omega^{T_{n'}}(du,dy)
\\&\leq&
\limsup_{n'\to\infty}\int_{\bbR^\sfp\times[0,1]}1_{\{u\leq q\}}y\>Q_\omega^{T_{n'}}(du,dy)
\\&\leq&
\int_{\bbR^\sfp\times[0,1]}1_{\{u\leq q\}}y\>Q_\omega(du,dy)
\eeas
for all $\omega\in\Omega_1$ and all $q\in\bbQ^\sfp$. 
If $I_\omega>0$, then 
\bea\label{170903-20}
\int_{\bbR^\sfp}1_{\{u<q\}}\widetilde{Q}_\omega(du)
&\leq&
\liminf_{n'\to\infty}\int_{\bbR^\sfp}1_{\{u\leq q\}}\widetilde{Q}_\omega^{T_{n'}}(du)
\nn\\&\leq&
\limsup_{n'\to\infty}\int_{\bbR^\sfp}1_{\{u\leq q\}}\widetilde{Q}_\omega^{T_{n'}}(du)
\nn\\&\leq&
\int_{\bbR^\sfp}1_{\{u\leq q\}}\widetilde{Q}_\omega(du)
\eea
for all $q\in\bbQ^\sfp$ and all $\omega\in\Omega_1$, where 
the probability measures $\widetilde{Q}_\omega^{T_{n'}}$ and $\widetilde{Q}_\omega$ on $\bbR^\sfp$ are given by 
\beas
\widetilde{{\sf Q}}(B) &=& 
\int_{\bbR^\sfp\times[0,1]}1_B(u)y\>{\sf Q}(du,dy)\bigg/\int_{\bbR^\sfp\times[0,1]}y\>{\sf Q}(du,dy)
\quad(B\in\bbB[\bbR^\sfp])
\eeas
for ${\sf Q}=Q_\omega^{T_{n'}}$ and $Q_\omega$. 
For any continuity point $r\in\bbR^\sfp$ of $\widetilde{Q}_\omega$, we take $q_1,q_2\in\bbR^\sfp$ with $q_1<r\leq q_2$ so that both are sufficiently 
close to $r$, and apply (\ref{170903-20}) to conclude  
$\widetilde{Q}_\omega^{T_{n'}}\to\widetilde{Q}_\omega$ for such $\omega$. Thus 
\bea\label{170903-21}
\int_{\bbR^\sfp\times[0,1]} f(u)y\>Q_\omega^{T_{n'}}(du,y)
&\to&
\int_{\bbR^\sfp\times[0,1]} f(u)y\>Q_\omega(du,y)
\eea
for $f\in C_b(\bbR^\sfp)$, $\omega\in\Omega_1$ with $I_\omega>0$. 
In the case $I_\omega=0$, it is obvious, so (\ref{170903-21}) holds for all $\omega\in\Omega_1$. 
This concludes the proof of Theorem \ref{170902-2}. 
\qed\halflineskip
\begin{en-text}
The already obtained (\ref{170902-3}) yields  
\bea\label{170903-7}
\int_{\hat{C}\times\bbR^\sfp\times[0,1]}F(x)y
1_{\{\Psi_T=1\}}1_{\{\Delta_q(x)\geq0\}}1_{\{u\leq q\}}
\>Q_T(dx,du,dy) &\to& 
\int_{\hat{C}\times\bbR^\sfp\times[0,1]}F(x)y\>Q(dx,du,dy)  
\eea
in probability  
for $F\in C_b(\hat{C})$, where 
$Q_T$ [resp. $Q$] is a regular conditional probability  of $(\tbbZ_T,\hat{u}_T,Y)$ [resp. $(\bbZ,\hat{u},Y)$] 
given $\calc$. 
Let ${\mathfrak F}$ be any countable subset of $\C_b(\hat{C})$. 
By subsequence argument, we can assume that the convergence (\ref{170903-7}) takes place 
for all $F\in{\mathfrak F}$ and 
for all $\omega\in\Omega_1\in\calf$ such that $P[\Omega_1]=1$, 
along a subsequence $(T_{n''})$ 
of any given sequence $(T_n)\nearrow\infty$. 
We may also assume that 
\beas 
\int_{\hat{C}\times\bbR^\sfp\times[0,1]}y\>Q_{T_{n''}}(dx,du,dy) &=& \int_{\hat{C}\times\bbR^\sfp\times[0,1]}y\>Q(dx,du,dy)
\eeas
for all $n''$ and $\omega\in\Omega_1$.

We will show 
\bea\label{17-903-8}
\int_{\hat{C}\times\bbR^\sfp\times[0,1]}f(u)y\>Q_{T_{n''}}(dx,du,dy) &\to& 
\int_{\hat{C}\times\bbR^\sfp\times[0,1]}f(u)y\>Q(dx,du,dy)
\eea
as $n''\to\infty$ for $\omega\in\Omega_1$. 
Let $\omega\in\Omega_1$. 
If $\int_{\hat{C}\times\bbR^\sfp\times[0,1]}y\>Q(dx,du,dy)=0$, then (\ref{17-903-8}) is obvious. 
Otherwise, define probability measures $\widetilde{Q}_{T_{n''}}$ and $\widetilde{Q}$ on $\hat{C}\times\bbR^\sfp$ by 
\beas
\widetilde{{\sf Q}}(B) &=& 
\int_{\hat{C}\times\bbR^\sfp\times[0,1]}1_B(x,u)y\>{\sf Q}(dx,du,dy)\bigg/\int_{\hat{C}\times\bbR^\sfp\times[0,1]}y\>{\sf Q}(dx,du,dy)
\quad(B\in\bbB[\hat{C}\times\bbR^\sfp])
\eeas
for ${\sf Q}=Q_{T_{n''}}$ and $Q$. 

For (\ref{170902-2}), it is sufficient to show 
\bea\label{17-903-6}
\int_{\hat{C}\times\bbR^\sfp\times\bbR}f(u)y\>Q_T(dx,du,dy) &\to^P& 
\int_{\hat{C}\times\bbR^\sfp\times\bbR}f(u)y\>Q(dx,du,dy)
\eea

We will show 
\bea
\int_{\bbR^{\sfp}\times I}
\eea

We may assume $Y\geq0$, furthermore $\int_{\hat{C}\times [0,1]}y\>dP^{T_{n'}}_\omega(dx,dy)>0$; 
otherwise . 
Set 
\beas 
\dot{P}_\omega^{T_{n'}}(B)=\int_{B\times [0,1]}y\>dP^{T_{n'}}_\omega(dx,dy)\big/\int_{\hat{C}\times [0,1]}y\>dP^{T_{n'}}_\omega(dx,dy) &\text{and}&
\dot{P}_\omega=P_\omega[\>\cdot\> Y]/E_\calc[Y]. 
\eeas
\koko

We will show 
\bea\label{17-903-5} 
E\big[f(\hat{u}_T)|\sigma[Y]\vee\calc\big] &\to^P& E\big[f(\hat{u})|\sigma[Y]\vee\calc\big]. 
\eea
For (\ref{17-903-5}), it is sufficient to show 
\bea\label{17-903-6}
\int_{\hat{C}\times\bbR^\sfp\times\bbR}f(u)yQ_T(dx,du,dy) &\to^P& 
\int_{\hat{C}\times\bbR^\sfp\times\bbR}f(u)yQ(dx,du,dy)
\eea
where $Q_T$ [resp. $Q$] is a regular conditional probability  of $(\tbbZ_T,\hat{u}_T,Y)$ [resp. $(\bbZ,\hat{u},Y)$] 
given $\sigma[Y]\vee\calc$, 

Therefore 
\beas
E_\calc\big[f(\hat{u}_T)Y\big] 
&=&
E_\calc\big[E_\calc\big[f(\hat{u}_T)|Y\big]\big] 
\>\to^P\> 
E_\calc\big[E_\calc\big[f(\hat{u})|Y\big]\big] 
\yeq 
E_\calc\big[f(\hat{u})Y\big]. 
\eeas
\end{en-text}

Conditional type PLD provides convergence of the conditional moments of $\hat{u}_T$ under truncation. 

\begin{theorem}\label{170904-2} %Let $L>\sfp$. 
Suppose that $q>0$ and $L>\sfp\vee q\vee1$. 
Suppose that $[A1]$, $[A2]$, $[A3]$ and $[A4]$ are satisfied. Then 
(\ref{170902-1.1}) of Theorem \ref{170902-1} holds. 
Moreover 
\bea\label{170904-1}
E_\calc\big[f(\hat{u}_T)Y\big]&\to^P&\overline{E}_\calc\big[f(\hat{u})Y\big]\quad 
\eea
for any bounded $\calg$-measurable random variable $Y$ and any $f\in C(\bbR^\sfp)$ such that 
$\limsup_{|u|\to\infty}|f(u)||u|^{-q}<\infty$. 
In particular, $\hat{u}_T\to^{d_{s}(\calg)}\hat{u}$ as $T\to\infty$. 
\end{theorem}
\begin{remark}\rm
The conditional expectation $\overline{E}_\calc[W]$ of a random variable $W$ is defined as the integral 
$\int_\bbR w\overline{P}^W_\omega(dw)$ with respect to a regular conditional probability $\overline{P}^W_\omega$ of $W$ given $\calc$. 
If $W\in L^1(\overline{P})$, then it coincides with the ordinary conditional expectation $\overline{E}[W|\calc]$ almost surely. However in general 
we do not assume $W\in L^1(\overline{P})$ nor $\overline{E}_\calc[W]\in L^1(\overline{P})$ in this article. 
We should be careful when applying the formula $\overline{E}_\calc[W\Psi_T]=\overline{E}_\calc[W]\Psi_T$; 
it is possible only when $\overline{E}_\calc[W]$ is well defined. 
The same remark applies to $E_\calc[W]$. 
On the other hand, each $\hat{u}_T$ is bounded because $\Theta$ is bounded, 
so $E_\calc[f(\hat{u}_T)Y]$ in (\ref{170904-1}) is well defined 
in any sense. 
%The truncation by $\Psi_T$ is necessary in (\ref{170904-1}). 
%Without it, the conditional integrability of the $f(\hat{u}_T)$ is not ensured. 
\end{remark}
\noindent
{\it Proof of Theorem \ref{170904-2}.} 
Let $p\in(q\vee1,L)$. We may assume $T\to\infty$ along $\calt\in{\mathfrak T}$. 
Almost surely 
\beas 
E_\calc\big[|\hat{u}_T|^p\Psi_T\big] 
&=&
\int_0^\infty pt^{p-1}E_\calc\big[1_{\{|\hat{u}_T|>t\}}\Psi_T\big]dt
\\&\leq&
1+p\sum_{r=1}^\infty(r+1)^{p-1}E_\calc\big[1_{\{|\hat{u}_T|>r\}}\Psi_T\big]
\\&\leq& 
1+p\sum_{r=1}^\infty(r+1)^{p-1}P_\calc\bigg[\sup_{u\in\bbV_T(r)}\bbZ_T(u)\geq1\bigg]\Psi_T
\\&\leq& 
1+p\sum_{r=1}^\infty(r+1)^{p-1}\bigg(1_{\{\ep(r)>1\}}+\frac{V_L}{r^L}\bigg)\>=:\>V
\eeas
where $V_L$ is a random variable bounding the right-hand side of the inequality of 
$[A1]$. 
The variable $V<\infty$ a.s. because $L>p$.  
By the convergence (\ref{170902-2}) of Theorem \ref{170903-30}, 
we have 
\beas
\overline{E}_\calc\big[|\hat{u}|^p\wedge A\big] 
&=& 
\lim_{n\to\infty}E_\calc\big[|\hat{u}_{T_n}|^p\wedge A\big] \Psi_{T_n}
\>\leq\> 
V\quad a.s.
\eeas
for $A\in\bbN$ and some sequence $(T_n)_{n\in\bbN}\up\infty$, 
and then the conditional monotone convergence theorem gives 
\beas
\overline{E}_\calc\big[|\hat{u}|^p\big] \leq V\quad a.s. 
\eeas
by letting $A\uparrow\infty$. 
For some constant $C>0$, $|f(u)|\leq C(1+|u|^q)$ for all $u\in\bbR^\sfp$. 
Let $f_A(u)=(-A)\vee f(u) \wedge A$ for $u\in\bbR^\sfp$ and $A>0$. 
Then for $\ep>0$, 
\beas &&
P\bigg[\bigg|E_\calc\big[f(\hat{u}_T)\Psi_TY\big]-\overline{E}_\calc\big[f(\hat{u})\Psi_TY\big]\bigg|>\ep\bigg]
\\&\leq&
P\bigg[\bigg|E_\calc\big[f_A(\hat{u}_T)Y\big]-\overline{E}_\calc\big[(f_A)(\hat{u})Y\big]\bigg|>\frac{\ep}{3}\bigg]
\\&&
+P\bigg[\bigg|E_\calc\big[|f-f_A|(\hat{u}_T)\Psi_T\big]\bigg|>\frac{\ep}{3(\|Y\|_\infty+1)}\bigg]
\\&&
+P\bigg[\bigg|\overline{E}_\calc\big[|f-f_A|(\hat{u})\big]\bigg|>\frac{\ep}{3(\|Y\|_\infty+1)}\bigg]. 
\eeas
We have 
\beas 
|f(u)-f_A(u)| 
&\leq& 
|f(u)|1_{\{|f(u)|\geq A\}}
\yleq
C(1+|u|^q)1_{\{C(1+|u|^q)\geq A\}}
\\&\leq&
C(1+|u|^p)\delta(A)
\eeas
for $A>C$ and  
\beas 
\delta(A) &=& \sup_{u:|u|\geq(A/C-1)^{1/q}}\frac{1+|u|^q}{1+|u|^p}.
\eeas
Then 
\beas 
\limsup_{T\to\infty}
P\bigg[\bigg|E_\calc\big[f(\hat{u}_T)\Psi_TY\big]-\overline{E}_\calc\big[f(\hat{u})\Psi_TY\big]\bigg|>\ep\bigg]
&\leq&
2P\bigg[C(1+V)\delta(A)>\frac{\ep}{3(\|Y\|_\infty+1)}\bigg].
\eeas
Since $\lim_{A\to\infty}\delta(A)=0$ and $\Psi_T\to^P1$, we obtain the convergence (\ref{170904-1}). 
\qed\halflineskip

\subsection{Quasi Bayesian estimator}

\bd
\im[[A1$^\sharp$\!\!]] %Conditional PLD under $\Psi_T\cdot P_\calc$ 
There exists $\calc$-measurable random variables $U$ and $V$ such that  
\beas 
P_\calc\bigg[S_T\bigg(r,\frac{U}{r^D}\bigg)\bigg]\Psi_T&\leq& \frac{V}{r^L}\quad a.s.\quad(T\in\calt,\>r\in\bbN)
\eeas
{\colred for every $\calt\in{\mathfrak T}$.} 
\ed

\begin{proposition}\label{170905-1} 
Suppose that $q\geq0$, $D>\sfp+q$ and $L>1$. 
Suppose that $[A1^\sharp]$, $[A2]$ and $[A3]$ are fulfilled. 
Then
\bea\label{170905-7}
E_\calc\bigg[f\bigg(\int_{\bbU_{T}}h(u)\bbZ_{T}(u)du\bigg)\>Y\bigg]
&\to^P&
E_\calc\bigg[f\bigg(\int_{\bbR^\sfp}h(u)\bbZ(u)du\bigg)\>Y\bigg]
\eea
as ${T}\to\infty$ 
for any $f\in C_b(\bbR^\sfk)$, any bounded $\calg$-measurable variable $Y$,  
and any $\bbR^\sfk$-valued measurable mapping $h$ satisfying $|h(u)|\leq C(1+|u|^q)$ for some constant $C$. 
\end{proposition}
\proof 
We may show the convergence along every sequence $\calt=(T_n)_{n\in\bbN}$ in 
${\mathfrak T}$. 
Choosing a sufficiently small positive constant $c_1$ depending on $(\sfp,D-q-\sfp)$, we obtain 
\bea\label{170905-5}&&
P_\calc\bigg[\int_{\{|u|\geq N\}\cap\bbU_T}|u|^q\bbZ_T(u)du>c_1UN^{-(D-q-\sfp)}\bigg]\Psi_T
\nn\\&\leq&
\sum_{r=N}^\infty P_\calc\bigg[{\colred r}%(r+1)
^{q+\sfp-1}
\sup_{u\in\bbV_T(r)}\bbZ_T(u)>Ur^{-(D-q-\sfp+1)}\bigg]\Psi_T
\nn\\&\leq&
\sum_{r=N}^\infty \frac{V}{r^L}\yleq(N-1)^{-(L-1)}V\quad a.s.
\eea
for $T>1$ and $N\in\bbN$.

For any $[0,1]$-valued $\calc$-measurable random variable $\Phi$ and sufficiently large $T$, 
\beas 
E\left[1_{\bigg\{\int_{\{|u|\geq N\}\cap B(0,R)}|u|^q\bbZ_T(u)du>c_1UN^{-(D-q-\sfp)}\bigg\}}\Psi_T\Phi\right]
&\leq&
(N-1)^{-(L-1)}E\big[V\Phi\big].
\eeas
Letting $T\to\infty$ with Theorem \ref{170902-1}, next letting $R\up\infty$, we obtain 
\bea\label{170905-6} 
P_\calc\bigg[\int_{\{|u|\geq N\}}|u|^q\bbZ(u)du>c_1UN^{-(D-q-\sfp)}\bigg]
&\leq&
(N-1)^{-(L-1)}V\quad a.s.
\eea

\begin{en-text}
In particular, 
\bea\label{170905-10} 
P_\calc\bigg[\int_{\bbU_T}|u|^q\bbZ_Tdu>c_1U+2^q|B(0,2)|\exp\big(kw_T(k^{-1},2)\big)\bigg]\Psi_T
&\leq& 
V\quad a.s.
\eea
for any $k\in\bbN$. 
\end{en-text}
We have 
\bea\label{170905-10}
\int_{\{|u|\leq N\}\cap\bbU_T}|u|^q\bbZ_T(u)du
&\leq&
N^q|B(0,N)|\exp\big(Nw_T(1,N)\big). 
\eea
In particular, this property is transferred to the limit as 
\bea\label{291008-5} 
\int_{\{|u|\leq N\}}|u|^q\bbZ(u) du
&\leq&
N^q|B(0,N)|\exp\big(Nw_T(1,N)\big)\quad a.s. 
\eea

Let $\ep>0$. 
Then by (\ref{170905-5}), (\ref{170905-10}), 
(\ref{170905-6}) and (\ref{291008-5}), 
there exists a number $K_0$ such that 
\beas
P\bigg[\sup_n P_\calc\bigg[\int_{\bbU_{T_n}}|h(u)|\bbZ_{T_n}(u)du>K_0\bigg]
<\frac{\ep}{4(1+\|f\|_\infty)}\bigg]
&\geq&
1-\ep
\eeas
and
\beas 
P\bigg[P_\calc\bigg[\int_{\bbU_{\bbR^\sfp}}|h(u)|\bbZ(u)du>K_0\bigg]%\Psi_{T_n}
\yleq
\frac{\ep}{4(1+\|f\|_\infty)}\bigg]
&>&
1-\ep
\eeas

Take $\eta\in(0,1)$ such that 
$|f(x_2)-f(x_1)|\leq \ep/2 $ for all $x_1,x_2\in\bbR^\sfk$ such that $|x_1|,|x_2|\leq K_0+1$ and $|x_2-x_1|\leq \eta$. 
%We may show the convergence along every sequence $(T_n)_{n\in\bbN}$ in $(0,\infty)$. 
{From} (\ref{170905-5}) and from (\ref{170905-6}), 
there exists $N_0\in\bbN$ such that 
\beas
P\bigg[\sup_{n}P_\calc\bigg[\int_{\{|u|\geq N_0\}\cap\bbU_{T_n}}|h(u)|\bbZ_{T_n}(u)du>\eta\bigg]%\Psi_{T_n}
\yleq
\frac{\ep}{4(1+\|f\|_\infty)}\bigg]
&>&
1-\ep
\eeas
and 
\beas
P\bigg[P_\calc\bigg[\int_{\{|u|\geq N_0\}\cap\bbR^\sfp}|h(u)|\bbZ(u)du>\eta\bigg]
\yleq
\frac{\ep}{4(1+\|f\|_\infty)}\bigg]
&>&
1-\ep
\eeas

\begin{en-text}
By [A2], there exists a $\delta_0>0$ such that 
\beas 
P\bigg[\sup_{n}P_\calc\big[W_{T_n}(\delta_0,N_0,\eta)\big]<\eta\bigg] &>& 1-\ep
\eeas
\end{en-text}
\begin{en-text}
Then there exists a constant $K_1=K_1(\ep)$ such that 
\beas
P\bigg[\sup_{n}P_\calc\bigg[\int_{\bbU_{T_n}}|h(u)|\bbZ_{T_n}(u)du>K_1\bigg]%\Psi_{T_n}
\yleq
\frac{\ep}{3(1+\|f\|_\infty)}\bigg]
&>&
1-\ep
\eeas
and that
\beas 
P\bigg[P_\calc\bigg[\int_{\bbU_{\bbR^\sfp}}|h(u)|\bbZ(u)du>K_1\bigg]%\Psi_{T_n}
\yleq
\frac{\ep}{3(1+\|f\|_\infty)}\bigg]
&>&
1-\ep
\eeas
\end{en-text}
Write 
\beas
J_T(S) &=& \int_{S\cap\bbU_T}h(u)\bbZ_T(u)du
\eeas
%for $S\subset\bbU_T$ 
and 
\beas
J(S) &=& \int_{S}h(u)\bbZ(u)du
\eeas
for $S\subset\bbR^\sfp$. 
Let %$A_0=A(\delta_0,N_0,\eta)$ and 
$B_0=B(0,N_0)$. We may assume $\|Y\|_\infty\leq1$. 
We will consider $n$ such that $B_0\subset\bbU_{T_n}$. 
Then
\beas&&
P\bigg[\bigg|E_\calc\big[f(J_{T_n}(\bbU_{T_n}))Y\big]-E_\calc\big[f(J_{T_n}(B_0))Y\big]\bigg|>\ep\bigg]
\\&\leq&
P\bigg[2\|f\|_\infty P_\calc\big[\big|J_{T_n}(\bbU_{T_n}-B_0)\big|>\eta\big]>\ep/2\bigg]
\\&&
+P\bigg[2\|f\|_\infty P_\calc\bigg[\int_{\bbU_{T_n}}\big|h(u)\big|\bbZ_{T_n}(u)du
>K_0\bigg]>\ep/2\bigg]
\\&<&
2\ep
\eeas
Similarly, 
\beas
P\bigg[\bigg|E_\calc\big[f(J(\bbR^\sfp))Y\big]-E_\calc\big[f(J(B_0))Y\big]\bigg|>\ep\bigg]
&<&
2\ep
\eeas
Now we apply Theorem \ref{170902-1} to the functional 
\beas 
F(x)&=& f\bigg(\int_{B_0}h(u)x(u)du\bigg)
\eeas
to obtain (\ref{170905-7}). 
\qed\halflineskip

\begin{en-text}
We have 
\beas &&
P\big[\big|E_\calc[J_T(\bbU_T)\Psi_TY]-\overline{E}_\calc[J(\bbR^\sfp)\Psi_TY]\big|>\ep\big]
\\&\leq&
P\big[\big|E_\calc[J_T(B(0,R))\Psi_TY]-\overline{E}_\calc[J(B(0,R))\Psi_TY]\big|>\ep/3\big]
\\&&
+P\big[\big|E_\calc[J_T(\bbU_T)Y]-E_\calc[J_T(B(0,R))Y]\big|>\ep/3\big]
\\&&
+P\big[\big|\overline{E}_\calc[J(B(0,R))Y]-\overline{E}_\calc[J(\bbR^\sfp)Y]\big|>\ep/3\big]
\eeas

uniform continuous on $|u|\leq A$

any $f\in C(\bbR^\sfp;\bbR^\sfk)$ satisfying $\sup_{u\in\bbR^\sfp}\big(1+|u|^q\big)^{-1}|h(u)|<\infty$. 
Let $\ep>0$. 

Since $|f(u)|\leq C(1+|u|^q)$ for some constant $C$, 
\beas
E_\calc\big[|J_T(S)|\Psi_T]
&\leq&
C\|\varpi\|_\infty\int_{S\cap\bbU_T}(1+|u|^q)E_\calc\big[\bbZ_T(u)\Psi_T\big]du
\\&\leq&
C(\sfp)\|\varpi\|_\infty
\bigg[e^{W_T(1,1)}+\sum_{r=1}^\infty\int_{S\cap\bbU_T\cap \{r\leq|u|<r+1\}}(1+|u|^q)
E_\calc\big[\bbZ_T(u)\Psi_T\big]du\bigg]
\eeas
\end{en-text}

Let $\tilde{u}_T=a_T^{-1}\big(\tilde{\theta}_T-\theta^*\big)$. 

\begin{theorem}\label{170904-3}
Let $D>\sfp+1$ and $L>1$. Suppose that $[A1^\sharp]$, $[A2]$ and $[A3]$ %and $[A5]$ 
are fulfilled. 
Then 
\beas 
E_\calc\big[f(\tilde{u}_T)Y\big] &\to^P& E_\calc\big[f(\tilde{u})Y\big] 
\eeas
for any $f\in C_b(\bbR^\sfp)$ and any bounded $\calg$-mesurable variable $Y$. 
In particular, $\tilde{u}_T\to^{d_{s}(\calg)}\tilde{u}$ as $T\to\infty$. 
\end{theorem}
\proof 
Let $J_T=\int_{\bbU_T}u\bbZ_T(u)\varpi(\theta^\dagger_T(u))du$, 
$I_T=\int_{\bbU_T}\bbZ_T(u)\varpi(\theta^\dagger_T(u))du$, 
$J_\infty=\int_{\bbR^\sfp}u\bbZ(u)\varpi(\theta^*)du$ and 
$I_\infty=\int_{\bbR^\sfp}\bbZ(u)\varpi(\theta^*)du$.  
Proposition \ref{170905-1} provides convergence 
${\colred (I_T,J_T)\to^{d_s}(I_\infty,J_\infty)}$. 
Therefore, for $\ep\in(0,1]$, we can find $a>0$ such that 
\beas
\limsup_{T\to\infty}
P\big[P_\calc\big[I_T\leq a\big]>
{\colred 3^{-1}(1+2\|f\|_\infty\|Y\|_\infty)^{-1}\ep}
%\ep/3
\big]
&\leq&
3\limsup_{T\to\infty}\ep^{-1}P\big[I_T\leq a\big]
\\&\leq&
3\ep^{-1}P\big[I_\infty\leq a\big]
\\&<& \ep/2.
\eeas
We have 
\beas
\bigg|E_\calc\big[f(\tilde{u}_T)Y\big]
-
E_\calc\left[ f\bigg(\frac{J_T}{I_T\vee a}\bigg)Y\right]\bigg|
&\leq&
2\|f\|_\infty \|Y\|_\infty P_\calc\left[I_T\leq a\right]\quad a.s.
\eeas
for $T\in(0,\infty]$, writng $\tilde{u}_\infty=\tilde{u}$.

Apply Proposition \ref{170905-1} to $q=1$ and the function 
$(x,y)\mapsto f(x/(y\vee a))$ ($x,y\in\bbR$) and $h(u)=(u,1)$, 
we obtain 
\beas
\limsup_{T\to\infty}P\bigg[\bigg|E_\calc\big[f(\tilde{u}_T)Y\big]-
E_\calc\big[f(\tilde{u})Y\big]\bigg|>\ep\bigg]
&<&
\ep. \qed
\eeas
\halflineskip

\bd
\im[[A5\!\!]] %Conditional PLD under $\Psi_T\cdot P_\calc$ 
%non-degeneracy
$\ds 
E_\calc\bigg[\frac{1}{\int_{\bbU_T}\bbZ_T(u)du}\bigg]
\yeq O_p(1)
$ 
as $T\to\infty$. 
\ed

\begin{theorem}\label{170904-4}
Suppose that $D>\sfp+q$, $q\geq1$ and $L>q+1$ and 
that $[A1^\sharp]$, $[A2]$, $[A3]$ and $[A5]$ 
are fulfilled. 
Then 
\beas
E_\calc\big[f(\tilde{u}_T)Y\big] &\to^P& 
\overline{E}_\calc\big[f(\tilde{u})Y\big]
\eeas
as $T\to\infty$ 
for any $\calg$-measurable bounded random variable $Y$ and any 
$f\in C(\bbR^\sfp;\bbR^\sfk)$ satisfying $\sup_{u\in\bbR^\sfp}\big(1+|u|^q\big)^{-1}|f(u)|<\infty$. 
\end{theorem}
\proof 
We may assume $\|Y\|_\infty\leq1$. 
Let $b\in(q,(D-\sfp)\wedge(L-1))$. 
{\colred On $\{\Psi_T=1\}$, }
\bea\label{291008-8}
E_\calc\big[|\tilde{u}_T|^bY\big] 
\nn&\leq& 
\sum_{r=0}^\infty E_\calc\bigg[
\frac{\int_{\{r<|u|\leq r+1\}\cap\bbU_T}|u|^b\bbZ_T(u)\varpi(\theta^\dagger_T(u))du}
{\int_{\bbU_T}\bbZ_T(u)\varpi(\theta^\dagger_T(u))du}
\bigg] 
\nn\\&\leq&
1+\sum_{r=1}^\infty(r+1)^b
\bigg\{P_\calc\bigg[\sup_{u\in\bbV_T(r)}\bbZ_T(u)\geq \frac{U}{r^D}\bigg]
+
C(\sfp,\varpi) U r^{\sfp-1-D}
E_\calc\bigg[\frac{1}{\int_{\bbU_T}\bbZ_T(u)du}\bigg]
\bigg\}
\nn\\&\leq&
1+C\bigg\{V\sum_{r=1}^\infty r^{b-L}
%+r^\sfp
+U\sum_{r=1}^\infty r^{\sfp+b-D-1}
E_\calc\bigg[\frac{1}{\int_{\bbU_T}\bbZ_T(u)du}\bigg]\bigg\}
\nn\\&=&
O_p(1)
\eea
as $T\to\infty$ under [A5]. 
Let $f_A=(-A)\vee f\wedge A$ for $A>0$. 
From (\ref{291008-8}), for any $\calt=(T_n)\in{\mathfrak T}$, 
\begin{en-text}
\beas 
P\bigg[\big|(f-f_A)(\tilde{u}_{T_n})\big|>\ep\bigg]
&\leq&
C\ep^{-q} A^{-(b-q)}E\bigg[
E_\calc (1+|\tilde{u}_T|^b)\bigg]. 
\eeas
\end{en-text}
\bea\label{291008-9} &&
\limsup_{n\to\infty}
E\bigg[\big|E_\calc\big[f(\tilde{u}_{T_n})Y\big]-E_\calc\big[f_A(\tilde{u}_{T_n})Y\big]\big|
\wedge 1\bigg]
\nn\\&\leq&
C\limsup_{n\to\infty}E\bigg[\big\{A^{-(b/q-1)}E_\calc\big[1+|\tilde{u}_{T_n}|^b\big]\big\}
\wedge 1\bigg]
\>\to\>0
\eea
as $A\to\infty$.
Moreover, by Theorem \ref{170904-3}, we have 
\bea\label{291008-10}
E_\calc\big[f_A(\tilde{u}_T)Y\big] &\to^P& 
\overline{E}_\calc\big[f_A(\tilde{u})Y\big]
\eea
as $T\to\infty$. 
Then (\ref{291008-9}) and (\ref{291008-10}) give the desired convergence. 
\qed\halflineskip

\begin{en-text}
%\\&\leq& 
%\infty\quad a.s.

\koko
\vspace{3cm}

Results: 

1. Convergence of $\bbZ_T\Psi_T$ under $P_\calc$ on $\bbR^\sfp$
\beas 
E_\calc\big[F(\bbZ_T)\Psi_T\big] &\to^{P}& \overline{E}_\calc\big[F(\bbZ)\big]%\quad a.s.
\eeas
Rewrite $E_\calc$ by a regular conditional probability on $\hat{C}$. Apply former QLA theory for any subsequence $(T_n)$. [This is important.
]

2. limit theorems and moment convergence under $\Psi_T\cdot P_\calc$. 

3. Apply $E[\  ]$ for bounded $f$. Obtain limit theorems for QLA estimators. \y
\end{en-text}

\begin{remark}\rm 
Localization is essential. If the effect of a slow component to the fast component is unbounded, 
sophisticated construction of $\Psi_T$ is required and 
{\colred any way using $E_\calc[f(\hat{u}_T)]$} 
without localization for unbounded $f$ is banned in general. 
\end{remark}

\begin{remark}\rm 
Generalization to the multi-scaling case is straightforward though we only treated a single scaling $a_T$. 
\end{remark}

\begin{en-text}
\begin{remark}\rm 
a transition kernel from $(\Omega,\calc)$ to $\hat{C}$
\beas 
\Omega\ni\omega\mapsto K_\omega\in\calp(\hat{C})
\eeas
\beas 
E_\calc\big[F(\bbZ_T)\big] &\to& \int_{\hat{C}}F(z)K(dz)
\eeas
for every $F\in C_b(\hat{C})$. 
\end{remark}

\begin{remark}\rm 
$\calc$-conditionally $\calg$-stable convergence
\end{remark}

{\colorr More concretely, conditions by moments here . }
\end{en-text}

\section{Conditional polynomial type large deviation}\label{290924-1}
As seen in Section \ref{170906-10}, the polynomial type large deviation inequality 
under conditional probability plays an essential role in 
the partial quasi likelihood analysis. 
We present a set of conditions that induces 
a conditional polynomial type large deviation (CPLD) inequality 
though there are various versions of sufficient conditions as \cite{yoshida2011polynomial} in unconditional cases. 
%It is possible by rephrasing it in the conditional setting. 

Suppose that $\bbH_T$ is of class $C^3$. 
$\lambda_{min}(A)$ and $\lambda_{max}(A)$ denote the minimum and maximum eigenvalues of a symmetric matrix $A$, repsectively. 
Let $L>0$ and let $b_T=(\lambda_{min}(a_T^\star a_T))^{-1}$. 
We assume that $b_T^{-1}\leq\lambda_{max}\big(a_T^\star a_T\big) \leq C_1b_T^{-1}$ for some constant $C_1\in[1,\infty)$. 
Let $\alpha\in(0,1)$ and $\beta=\alpha/(1-\alpha)$. 
Let $\rho$ be a positive constant; practically $\rho=2$ in most cases. 

\bd 
\im[[B1\!\!]] 
Parameters $\beta_1$, $\rho_1$, $\rho_2$ and $\beta_2$ satisfy the following inequalities 
\beas &&
0<\beta_1<\half,\quad 0<\rho_1<\min\bigg\{1,\beta,\frac{2\beta_1}{1-\alpha}\bigg\},
\\&&
\alpha\rho<\rho_2,\quad\beta_2\geq0,\quad 1-2\beta_2-\rho_2>0.
\eeas
\ed
For example, the following two sets of conditions respectively satisfy [B1]. 
\bd\im[(i)] 
$\beta_1	=\alpha/2$, $\rho_1=\alpha$, $\rho_2=3\alpha$ and  $\beta_2=\alpha$ for 
$\rho=2$ and $\alpha\in(0,1/5)$. 
\im[(ii)]
$\beta_1	=\alpha/2$, $\rho_1=\alpha$, $\rho_2=3\alpha$ and  $\beta_2=0$ for 
$\rho=2$ and $\alpha\in(0,1/3)$. 
\ed

Let $\bbY:\Omega\times\Theta\to\bbR$ be a random field, i.e., a measurable mapping. 
\bd 
\im[[B2\!\!]] There exists a positive random variable $\chi_0$ satisfying the following conditions. 
\bd\im[(i)] 
With probability one, 
\beas 
\bbY(\theta)\yeq\bbY(\theta)-\bbY(\theta^*)&\leq&-\chi_0\big|\theta-\theta^*\big|^\rho
\eeas
for all $\theta\in\Theta$. 
\im[(ii)] 
$\ds  
\sup_{r\in\bbN}r^LP_\calc\big[\chi_0\leq r^{-(\rho_2-\alpha\rho)}\big] \><\> \infty$ 
a.s.
%For some $\calc$-measurable real random variable $V_{L,2}$, 
%\beas 
%P_\calc\big[\chi_0\leq r^{-(\rho_2-\alpha\rho)}\big] &\leq& \frac{V_{L,2}}{r^L}\quad a.s.
%\eeas
%for all $r>1$. 
\ed
\ed

Let $\Gamma$ be a $\sfp\times\sfp$ positive definite random matrix. 
%$\lambda_{min}(\Gamma)$ denotes the minimum eigenvalue of $\Gamma$. 
%
\bd
\im[[B3\!\!]] 
$\sup_{r\in\bbN}r^LP_\calc\big[\lambda_{min}(\Gamma)< 4r^{-\rho}\big]
\><\>\infty$ a.s.
\ed
%There exists a $\calc$-measurable real random variable $V_{L,3}$ such that 
%\beas 
%P_\calc\bigg[\lambda_{min}\Gamma< 4r^{-\rho}\bigg]
%&\leq&\frac{V_{L,3}}{r^L}
%\eeas
%\ed
\halflineskip

Let 
\beas 
\bbY_T(\theta) &=& \frac{1}{b_T}\big(\bbH_T(\theta)-\bbH_T(\theta^*)\big).
\eeas
Define a $\sfp$-dimensional random variable $\Delta_T$ and a $\sfp\times\sfp$ random matrix $\Gamma_T$ by 
\beas 
\Delta_T[u] &=& \partial_\theta\bbH_T(\theta^*)[a_Tu]\quad(u\in\bbR^\sfp)
\eeas
and 
\beas 
\Gamma_T[u^{\otimes2}] &=& -\partial_\theta^2\bbH_T(\theta^*)\big[\big(a_Tu\big)^{\otimes2}\big]\quad(u\in\bbR^\sfp)
\eeas
respectively.

\bd 
\im[[B4\!\!]] 
For $M_1=L(1-\rho_1)^{-1}$, 
\beas
\sup_{T\in \calt}E_\calc\big[\big|\Delta_T\big|^{M_1}\big]{\colorr \Psi_T}&<&\infty\quad a.s.
\eeas
for every $\calt\in{\mathfrak T}$. Moreover, for $M_2=L(1-2\beta_2-\rho_2)^{-1}$, 
\beas 
\sup_{T\in \calt}E_\calc\bigg[\bigg(\sup_{\theta\in\Theta}b_T^{\half-\beta_2}\big|\bbY_T(\theta)-\bbY(\theta)\big|\bigg)^{M_2}\bigg]{\colorr \Psi_T}
&<&\infty\quad a.s.
\eeas
for every $\calt\in{\mathfrak T}$.
\ed

\bd 
\im[[B5\!\!]] 
For $M_3=L(\beta-\rho_1)^{-1}$, 
\beas 
\sup_{T\in \calt}E_\calc\bigg[\bigg(b_T^{-1}\sup_{\theta\in\Theta}
\big|\partial_\theta^3\bbH_T(\theta)\big|\bigg)^{M_3}\bigg]{\colorr \Psi_T}
&<&\infty\quad a.s.
\eeas
for every $\calt\in{\mathfrak T}$.
Moreover, for $M_4=L\big(2\beta_1(1-\alpha)^{-1}-\rho_1\big)^{-1}$, 
\beas 
\sup_{T\in \calt}E_\calc\bigg[\big(b_T^{\beta_1}\big|\Gamma_T-\Gamma\big|\big)^{M_4}\bigg]{\colorr \Psi_T}
&<&\infty\quad a.s.
\eeas
for every $\calt\in{\mathfrak T}$.
\ed

\begin{theorem}\label{290916-1}
Suppose that Conditions $[B1]$-$[B5]$ are satisfied. Then 
\bea\label{290916-3}
\sup_{T\in\calt}\sup_{r\in\bbN}r^LP_\calc\bigg[\sup_{u\in\bbV_T(r)}\bbZ_T(u)\geq \exp\big(-2^{-1}r^{2-(\rho_1\vee\rho_2)}\big)\bigg]\Psi_T
&<&
\infty\quad a.s.
\eea
for every $\calt\in{\mathfrak T}$. 
Moreover, $\bbZ_T$ has a LAQ representation 
\beas
\bbZ_T(u) &=&\exp\bigg(
\Delta_T[u]-\half \Gamma[u^{\otimes2}]+r_T(u)\bigg)
\eeas
with $r_T(u)\to^P0$ as $T\to\infty$ for every $u\in\bbR^\sfp$. 
\end{theorem}
\proof Arbitrarily given $\calt\in{\mathfrak T}$, we will follow the proof of Theorems 1 and 2 of Yoshida \cite{yoshida2011polynomial} 
under $P_\calc[\>\cdot\>\Psi_T]$. This is valid because positivity of the expectation used in the proof in \cite{yoshida2011polynomial} 
is obviously valid for $P_\calc[\>\cdot\>\Psi_T]$ a.s. in the present case. 
Condition $[A4']$ of \cite{yoshida2011polynomial} is the present Condition $[B1]$. 
Condition $[A1'']$ of \cite{yoshida2011polynomial} is satisfied by $[B5]$ in conditional version. 
Then as Lemma 1 of \cite{yoshida2011polynomial} deduced $[A1']$ therein, we obtain 
\bea\label{290917-1} 
\sup_{T\in\calt}\sup_{r\in\bbN}r^LP_\calc\big[S'_T(r)^c\big]
{\colred \Psi_{T}} &<& \infty\quad a.s.
\eea
where 
\beas
S'_T(r) &=& 
\left\{\omega;\>\sup_{\theta\in\Theta\atop
b_T^{-1/2}r\leq|\theta-\theta^*|\leq C_1^{1/2}\delta_T}
\big|\Gamma_T(\theta)-\Gamma\big|<\ep_1(r)\right\}
\eeas
with $\delta_T=b_T^{-\alpha/2}$ and $\ep_1(r)=r^{-\rho_1}$. 
$\Gamma_T(\theta)$ is defined by 
\beas 
\Gamma_T(\theta)[u^{\otimes2}] &=& -\partial_\theta^2\bbH_T(\theta)\big[\big(a_Tu\big)^{\otimes2}\big]\quad(u\in\bbR^\sfp).
\eeas

The variable $r_T(u)$ is defined by 
the LAQ representation of $\bbZ_T(u)$: 
\bea\label{290922-1} 
\bbZ_T(u) &=&\exp\bigg(
\Delta_T[u]-\half \Gamma[u^{\otimes2}]+r_T(u)\bigg)
\eea
where $r_T(u)$ admits the expression 
\beas 
r_T(u) &=& \int_0^1 (1-s)\bigg\{\Gamma[u^{\otimes2}]-\Gamma_T(\theta^\dagger_T(su))[u^{\otimes2}]\bigg\}ds
\eeas
for every $u\in\bbR^\sfp$ and sufficiently large $T$ depending on $u$. 
Then from (\ref{290917-1}), we obtain, as a counterpart of $[A1]$ of \cite{yoshida2011polynomial}, 
\bea\label{290917-2}
\sup_{T\in\calt}\sup_{r\in\bbN}r^LP_\calc\big[S_T(r)^c\big]{\colred \Psi_{T}}  &<& \infty\quad a.s.
\eea
where 
\beas
S_T(r) &=& \left\{\omega;\>\sup_{u\in \bbU_T(r)}(1+|u|^2)^{-1}\big|r_T(u)\big|<\ep_1(r)\right\}
\eeas
with 
$\bbU_T(r)=\big\{u\in\bbU_T;\>r\leq|u|\leq\delta_Tb_T^{1/2}\big\}$. 
%The event $S_T(r)$ is well defined (i.e. $r_T(u)$ is well defined) for large $T$. 

Condition $[B3]$ serves as $[A2]$ in \cite{yoshida2011polynomial}, 
Condition $[B2]$ (i) as $[A3]$, Condition $[B2]$ (ii) as $[A5]$, and 
Condition $[B4]$ as $[A6]$, respectively. 
As already mentioned, $[B1]$ is $[A4']$ that is stronger than $[A4]$. 
Now by using (\ref{290917-2}) and $[B1]$-$[B4]$, 
we follow the line of the proof of Theorem 1 of \cite{yoshida2011polynomial} given under $[A1]$-$[A6]$,  
to obtain (\ref{290916-3}). 

Condition $[B5]$ implies $r_T(u)\to^P0$ as $T\to\infty$. 
\qed\halflineskip

Obviously the inequality (\ref{290916-3}) ensures $[A1^\sharp]$, $[A1]$ and $[A1^\flat]$.

\section{Partial mixing}\label{290924-2}
Partial mixing is a structure we often meet in applications of the partial quasi likelihood analysis, though it is not the all. We state a Rosenthal type inequality 
under conditional expectation. 
\begin{lemma}\label{290824-1}
Let $2\leq p<r$. 
Given a probability space $(\Omega,\calf,P)$ 
and a sub $\sigma$-fields $\calc$ of $\calf$, 
let $\calg_j$ and $\calh_j$ $(j=1,2,...,n)$ 
be sub $\sigma$-fields of $\calf$ such that $\calg_j\cap\calh_j\supset\calc$ 
for all $j=1,...,n$. 
Let $X=(X_j)_{j=1,...,n}$ be a sequence of random variables 
such that $X_j\in L_r(\Omega,\calg_j\cap\calh_j,P)$ and $E_\calc[X_j]=0$ a.s. 
Suppose that $[0,1/2]$-valued $\calc$-measurable random variables $\alpha_\calc(h)$ 
satisfy 
\beas 
\alpha_\calc(h) &\geq&
\sup_{k=1,...,n-h}\sup\big\{\big|P_\calc[A\cap B]-P_\calc[A]P_\calc[B]\big|;\> A\in\calg_k,\>B\in\calh_{k+h}\big\}
\quad 
\eeas
for $h=1,...,n-1$ on some $\Omega_0\in\calf$ such that $P[\Omega_0]=1$. 
Then
\beas 
E_\calc\bigg[\max_{k=1,..,n}\bigg|\sum_{j=1}^kX_j\bigg|^p\bigg]
&\leq&
C(p,r)\max_{j=1,...,n}E_\calc\big[\big|X_j\big|^r\big]^{p/r}\>
\\&&
\times\bigg[n^{p/2}\bigg(1+\sum_{h=1}^{n-1}\alpha_\calc(h)^{1-2/r}\bigg)^{\colorr p/2}
+n\sum_{h=1}^{n-1}(h+1)^{p-2}\alpha_\calc(h)^{1-p/r}\bigg]\quad a.s.
\eeas
where $C(p,r)$ is a constant depending only on $p$ and $r$. 
\end{lemma}
\proof 
Denoted by $Q_j$ the random upper quantile function of 
a regular conditional distribution $P_\omega^{|X_j|}$ of $|X_j|$ given $\calc$, 
i.e., an inverse of the function $t\mapsto P_\omega^{|X_j|}[(t,\infty)]$. 
Let $\alpha_\calc(0)=1/2$. 
The random function $\alpha_\calc^{-1}:(0,1)\to\bbZ_+$ is defined by 
$\alpha_\calc^{-1}(u)=\sum_{h\in\bbZ_+}1_{\{u<\alpha_\calc(h)\}}$. 
We apply Theorem 6.3 of Rio \cite{rio2017asymptotic} under the conditional probability 
$P_\calc$ to obtain 
\bea\label{170827-1} 
E_\calc\bigg[\max_{k=1,..,n}\bigg|\sum_{j=1}^kX_j\bigg|^p\bigg]
&\leq&
C(p)\bigg[s_\calc(n)^p+n\int_0^1\big[\alpha_\calc^{-1}(u)\wedge n\big]^{p-1}Q_\calc(u)^pdu\bigg]
\quad a.s.
\eea
where $C(p)$ is a constant depending only on $p$,  
$Q_\calc=\max_{j=1,...,n}Q_i$ is %an $[0,\infty]$-valued 
a $\calc$-measurable random function of $u\in(0,1)$  
and 
\beas 
s_\calc(n)^2 &=& 
\sum_{i,j=1}^n\big|\text{Cov}_\calc[X_i,X_j]\big|,
\eeas
$\text{Cov}_\calc$ denoting $\calc$-conditional covariance. 
We shall estimate the right-hand side of (\ref{170827-1}). 

Since $Q_\calc(u)\leq \rho u^{-1/r}$ 
with $\rho=\max_{j=1,...,n}E_\calc[|X_j|^r]^{1/r}$, we have 
\beas 
\int_0^1 \big[\alpha_\calc^{-1}(u)\wedge n\big]^{p-1}Q_\calc(u)^pdu
&\leq&
\int_0^{1/2}
 \bigg[\bigg(\sum_{h\in\bbZ_+}1_{\{u<\alpha_\calc(h)\}}\bigg)\wedge n\bigg]^{p-1}
\rho^p u^{-p/r}du
\\&\leq&
\int_0^{\alpha_\calc(n)}n^{p-1}\rho^p u^{-p/r}du
+
\sum_{\ell=1}^n\int_{\alpha_\calc(\ell)}^{\alpha_\calc(\ell-1)}\ell^{p-1}\rho^p u^{-p/r}du
\\&\leq&
\rho^p\big(1-p/r\big)^{-1}\sum_{\ell=1}^n\ell^{p-2}\alpha_\calc(\ell-1)^{1-p/r}
\eeas
and move the term for $\ell=1$ into the error bound we will consider. 

By the covariance inequality applied to the conditional situation, 
\beas 
s_\calc(n)^2
&\leq&
\sum_{i,j=1}^n2\alpha_\calc(|i-j|)^{1-2/r}\rho^2
\>\leq\>
4\rho^2n\sum_{h=0}^{n-1}\alpha_\calc(h)^{1-2/r}.
\eeas
Bring the above two estimates into (\ref{170827-1}), we complete the proof. 
\begin{en-text}
We give a proof for convenience of reference. 
First, 
\beas 
E_\calc\bigg[\bigg|n^{-1}\sum_{j=1}^nX_j\bigg|^q1_{\cup_{j=1}^n\{|X_j|\geq n^\ep\}}\bigg]
&\leq&
E_\calc\bigg[\bigg(n^{-1}\sum_{j=1}^n|X_j|\bigg)^r\bigg]^{q/r}
P_\calc\bigg[\bigcup_{j=1}^n\{|X_j|\geq n^\ep\}\bigg]^{1-q/r}
\\&\leq&
E_\calc\bigg[n^{-1}\sum_{j=1}^n|X_j|^r\bigg]^{q/r}
\bigg(n^{1-r\ep}E_\calc\big[|X_j|^r\big]\bigg)^{1-q/r}
\\&\leq&
n^{(1-r\ep)(1-q/r)}\max_{j=1,...,n}E_\calc[|X_j|^r] 
\eeas
almost surely. 
Next, writing $1_j=1_{\{|X_j|< n^\ep\}}$, we have 
\beas 
E_\calc\bigg[\bigg|n^{-1}\sum_{j=1}^nX_j1_j\bigg|^q\bigg]
&\leq&
E_\calc\bigg[\bigg|n^{-1}\sum_{j=1}^nX_j1_j\bigg|^2\bigg]n^{(q-2)\ep}
\\&\leq&
n^{-2+(q-2)\ep}\sum_{i,j=1}^n\big\{\text{Cov}_\calc[X_i1_i,X_j1_j]+E_\calc[X_i1_i]E_\calc[X_j1_j]\big\}
\\&\leq&
n^{-2+(q-2)\ep}\sum_{i,j=1}^n\big\{2n^{2\ep}\alpha^X_\calc(|i-j|)+n^{-2r\ep}E_\calc[|X_i|^r]E_\calc[|X_j|^r]\big\}
\eeas
almost surely 
by the conditional version of Ibragimov's covariance inequality. 
Then the desired inequality follows from these inequalities. 
\end{en-text}
\qed\y

\begin{en-text}
\section{An intensity model}
\bi 
\im A Cox type model
\bea\label{20161220-1} 
\lambda^M(t,\theta^M) &=& 
\lambda_0(t)\exp\bigg(\sum_{j}\theta^M_jX^M_j(t)\bigg)
\>=\>\lambda_0(t)\bar{\lambda}^M(t,\theta^M),
\nn\\
\lambda^L(t,\theta^L) &=& 
\lambda_0(t)\exp\bigg(\sum_{k}\theta^L_kX^L_k(t)\bigg)
\>=\>\lambda_0(t)\bar{\lambda}^L(t,\theta^L).
\eea
We assume that a (possibly random) baseline hazard $\lambda_0(t)$ is common for the two intensities. 
Intercepts can be covariates as \sout{$X^M_0(t)=1$} 
{\colorr $X^M_0(t)=0$} and $X^L_0(t)=1$. 

\noindent
{\colorr Remark. If we took $X^M_0(t)=1$, then $(\theta^M_0,\theta^L_0)$ would be 
non-identifiable due to freedom of $\lambda_0(t)$ unless $ \lambda_0(t)$ was constrained. 
In this case, $\lambda_0(t)$ is the baseline intensity of the market order, and 
the symmetry of the modeling is broken. }

\im quasi log likelihood
\beas 
\ell^M_T(\theta^M)
&=&
\int_0^T \log\lambda^M(t,\theta^M)dN^M_t-\int_0^T\lambda^M(t,\theta)dt
\\&=&
\int_0^T \log\big\{\lambda_0(t)\bar{\lambda}^M(t,\theta^M)\big\}dN^M_t-\int_0^T\lambda_0(t)\bar{\lambda}^M(t,\theta^M)dt.
\eeas

\im score function
\beas 
\partial_{\theta^M}\ell^M_T(\theta^M)
&=&
\int_0^T \bar{\lambda}^M(t,\theta^M)^{-1}\partial_{\theta^M}\bar{\lambda}^M(t,\theta^M)dN^M_t-\int_0^T\lambda_0(t)\partial_{\theta^M}\bar{\lambda}^M(t,\theta^M)dt
\\&=&
\int_0^T \bar{\lambda}^M(t,\theta^M)^{-1}\partial_{\theta^M}\bar{\lambda}^M(t,\theta^M)
\bigg(dN^M_t-\lambda_0(t)\bar{\lambda}^M(t,\theta^M)dt\bigg).
\eeas

\im Estimating functions (E1)
\beas 
\bbF^M_T(\theta^M,\theta^L)
&=&
\int_0^T \bar{\lambda}^M(t,\theta^M)^{-1}\partial_{\theta^M}\bar{\lambda}^M(t,\theta^M)
\bigg(dN^M_t-\bar{\lambda}^L(t,\theta^L)^{-1}
\bar{\lambda}^M(t,\theta^M)dN^L_t\bigg)
\\&=&
\int_0^T \partial_{\theta^M}\bar{\lambda}^M(t,\theta^M)
\bigg(\bar{\lambda}^M(t,\theta^M)^{-1}dN^M_t-\bar{\lambda}^L(t,\theta^L)^{-1}
dN^L_t\bigg)
\>=\>0.
\eeas

\beas 
\bbF^L_T(\theta^M,\theta^L)
&=&
\int_0^T \partial_{\theta^L}\bar{\lambda}^L(t,\theta^L)
\bigg(\bar{\lambda}^L(t,\theta^L)^{-1}dN^L_t-\bar{\lambda}^M(t,\theta^M)^{-1}
dN^M_t\bigg)
\>=\>0.
\eeas

\noindent 
{\colorr Remark. The derivative of $(\bbF^M_T,\bbF^L_T)$ in $(\theta^M,\theta^L)$ 
is not symmetric although the elements are similar to those of the Fisher information matrix. 
Nondegeneracy of the derivative is not clear. 
The asymptotic variance of the estimator is of so-called sandwich type. 
The estimator (E1) as well as (E2) is quite ad hoc. 
More consideration is necessary. }

\im score function
\beas 
\partial_{\theta^M}\ell^M_T(\theta^M)
&=&
\int_0^T \bar{\lambda}^M(t,\theta^M)^{-1}\partial_{\theta^M}\bar{\lambda}^M(t,\theta^M)dN^M_t-\int_0^T\lambda_0(t)\partial_{\theta^M}\bar{\lambda}^M(t,\theta^M)dt
\\&=&
\int_0^T \bar{\lambda}^M(t,\theta^M)^{-1}\partial_{\theta^M}\bar{\lambda}^M(t,\theta^M)
\bigg(dN^M_t-\lambda_0(t)\bar{\lambda}^M(t,\theta^M)dt\bigg).
\eeas

\im Estimating functions (E2)
\beas 
\bbG^M_T(\theta^M,\theta^L)
&=&
\int_0^T \bar{\lambda}^M(t,\theta^M)^{-1}\partial_{\theta^M}\bar{\lambda}^M(t,\theta^M)
\bigg(dN^M_t-\{\bar{\lambda}^M(t,\theta^M)+\bar{\lambda}^L(t,\theta^L)\}^{-1}\bar{\lambda}^M(t,\theta^M)
\{dN^M_t+dN^L_t\}\bigg)
\\&=&
\int_0^T \bar{\lambda}^M(t,\theta^M)^{-1}\partial_{\theta^M}\bar{\lambda}^M(t,\theta^M)
\bigg(
\frac{\bar{\lambda}^L(t,\theta^L)}{\bar{\lambda}^M(t,\theta^M)+\bar{\lambda}^L(t,\theta^L)}
dN^M_t
-\frac{\bar{\lambda}^M(t,\theta^M)}{\bar{\lambda}^M(t,\theta^M)+\bar{\lambda}^L(t,\theta^L)}
dN^L_t\bigg)
\\&=&
\int_0^T \frac{\bar{\lambda}^M(t,\theta^M)^{-1}\partial_{\theta^M}\bar{\lambda}^M(t,\theta^M)}
{\bar{\lambda}^M(t,\theta^M)+\bar{\lambda}^L(t,\theta^L)}
\bigg(
\bar{\lambda}^L(t,\theta^L)dN^M_t
-\bar{\lambda}^M(t,\theta^M)
dN^L_t\bigg)
\>=\>0.
\eeas

\beas 
\bbG^L_T(\theta^M,\theta^L)
&=&
\int_0^T \frac{\bar{\lambda}^L(t,\theta^L)^{-1}\partial_{\theta^L}\bar{\lambda}^L(t,\theta^L)}
{\bar{\lambda}^M(t,\theta^M)+\bar{\lambda}^L(t,\theta^L)}
\bigg(
\bar{\lambda}^M(t,\theta^M)dN^L_t
-\bar{\lambda}^L(t,\theta^L)dN^M_t\bigg)
\>=\>0.
\eeas

\im Hybrid estimating functions
\bi
\im $(\bbF^M_T(\theta^M,\theta^L),\bbG^L_T(\theta^M,\theta^L))=0$
\im $(\bbF^L_T(\theta^M,\theta^L),\bbG^M_T(\theta^M,\theta^L))=0$
\ei

\im an estimator $\hat{\lambda}_0(t;T)$ for the baseline intensity $\lambda_0(t)$ by the data on $[0,T]$
\beas 
\hat{\lambda}_0(t;T) dt 
&=& 
\frac{rdN^M_t+(1-r)dN^L_t}{r\bar{\lambda}^M(t,\hat{\theta}^M_T)
+(1-r)\bar{\lambda}^L(t,\hat{\theta}^L_T)},
\eeas
where $r$ is a parameter, e.g., $r=1/2$.

\im $L$ can be the market order of other securities. 
Example. 
\beas 
\lambda^{i,M}(t,\theta^{i,M}) &=& 
\lambda_0(t)\exp\bigg(\sum_{j}\theta^{i,M}_jX^{i,M}_j(t)\bigg)
\>=\>\lambda_0(t)\bar{\lambda}^{i,M}(t,\theta^{i,M})
\quad(i=1,...,I).
\eeas
\bi
\im Estimating functions similar to (E2) can be constructed. 
\im Covariates that are common among $i$'s can be introduced. 
It is possible to consider common parameters $\theta^{i_1,M}_{j_1}=\theta^{i_2,M}_{j_2}=...$. 
\im Example. $(1,M)$ for the market order on the bid side, $(2,M)$ for the market order on the ask side. 
\ei

\im modeling $\lambda_0(t)$ 
\bi
\im switching model 
\im ?
\ei
\ei

\newpage
\section{Discussion}\label{20161220-2}

We are working with the model (\ref{20161220-1}): 
\beas
\l\{\begin{array}{rcl}
\lambda^M(t,\theta^M) &=& 
\lambda_0(t)\exp\bigg(\sum_{j}\theta^M_jX^M_j(t)\bigg)
\>=\>\lambda_0(t)\bar{\lambda}^M(t,\theta^M),
\nn\y
\lambda^L(t,\theta^L) &=& 
\lambda_0(t)\exp\bigg(\sum_{k}\theta^L_kX^L_k(t)\bigg)
\>=\>\lambda_0(t)\bar{\lambda}^L(t,\theta^L).
\end{array}\r.
\eeas
The processes $X^M_j(t)$ and $X^L_k(t)$ are assumed to be observable, on the other hand, 
$\lambda_0(t)$ is unobservable in general. 

The model 
\bea\label{20161222-1}
\l\{\begin{array}{rcl}
\lambda^M(t,\theta^M) &=& \lambda_0\exp\big(\theta^M_1X(t)\big),
\y
\lambda^L(t,\theta^L)&=&\lambda_0\exp\big(\theta^L_0+\theta^L_1X(t)\big)
\end{array}\r.
\eea
can be rewritten as e.g. 
\beas 
\l\{\begin{array}{rcl}
\lambda^M(t,\theta^M) &=& \lambda_0\exp\big(\theta^L_1X(t)\big)
\times \exp\big((\theta^M_1-\theta^L_1)X(t)\big),
\y
\lambda^L(t,\theta^L)&=&\lambda_0\exp\big(\theta^L_1X(t)\big)
\times\exp\big(\theta^L_0\big)
\end{array}\r.
\eeas
with 
$\lambda_0(t)=\lambda_0\exp\big(\theta^L_1X(t)\big)$. 
Therefore, for this model, it is natural to obtain consistent estimates 
only for  $\theta^M_1-\theta^L_1$ not for $(\theta^M_1,\theta^L_1)$. 
Indeed, the estimating functions 
\beas 
\l\{\begin{array}{rcl}
F^M_1(\theta^M,\theta^L) &=& \sum_{t^M_i}X(t^M_i-)
-\sum_{t^L_i}X(t^L_i-)\exp\big(\theta^M_1X(t^L_i-)-\theta^L_0-\theta^L_1X(t^L_i-)\big)
\y
F^L_0(\theta^M,\theta^L) &=& \sum_{t^L_i}1
-\sum_{t^M_i}1\cdot\exp\big(-\theta^M_1X(t^M_i-)+\theta^L_0+\theta^L_1X(t^M_i-)\big)
\y
F^L_1(\theta^M,\theta^L) &=& \sum_{t^L_i}X(t^L_i-)
-\sum_{t^M_i}X(t^M_i-)\exp\big(-\theta^M_1X(t^M_i-)+\theta^L_0+\theta^L_1X(t^M_i-)\big)
\end{array}\r.
\eeas
are functions of $\theta^L_0$ and $\theta^M_1-\theta^L_1$. 

In the present situation, since $X(t)$ is observable, 
once these parameters were consistently estimated, the parameters $\lambda_0$ and 
$\theta^L_1$ (and hence $\theta^M_1$) may be estimated by maximizing 
\beas 
\ell_T(\lambda_0,\theta^L_1) 
&=& 
\int_0^T\log \lambda^L(t,\hat{\theta}^L_0,\theta^L_1)dN^L_t
-\int_0^T \lambda^L(t,\hat{\theta}^L_0,\theta^L_1)dt,
\eeas
where $\hat{\theta}^L_0$ is the estimator obtained above based on the data 
on $[0,T]$. 

As a matter of fact, the quasi likelihood approach can be taken for the model (\ref{20161222-1}). 

\end{en-text}
In the following two sections, we will present applications of 
the partial quasi likelihood analysis. 

\section{Diffusion process having a component with a slow mixing rate %long range dependence
}\label{290924-5}
\subsection{Partial QLA for a stochastic regression model}
Given a stochastic basis $(\Omega,\calf,{\bf F},P)$, 
${\bf F}=(\calf_t)_{t\in\bbR_+}$, we consider a stochastic regression model 
\bea\label{290920-1}
Y_t &=& Y_0+\int_0^t b(X_s,\theta^*)ds 
+ \int_0^t \sigma(X_s)dw_s\quad(t\in\bbR_+).
\eea
Here  
%$\gamma$ expresses a measurable random scenery taking values in a measurable space $({\sf G},\bbB_{\sf G})$. 
$X=(X_t)_{t\in\bbR_+}$ is a stochastic process taking values in a measurable space $({\sf X},\bbB_{\sf X})$. 
$\Theta$ is a bounded domain in $\bbR^\sfp$. 
We assume that the boundary of $\Theta$ is as good as it admits the ordinary Sobolev's inequality for 
the embedding $W^{1,p}(\Theta)\inclusion C_b(\Theta)$ 
for $p>\sfp$. 
Moreover, 
$b:{\sf X}\times\Theta\to\bbR^\sfm$ and 
$\sigma:{\sf X}\to \bbR^\sfm\otimes\bbR^\sfr$ 
are given functions. 
$w=(w_t)_{t\in\bbR_+}$ is an $\sfr$-dimensional standard 
${\bf F}$-Wiener process. 
We assume that $b(X_t,\theta)$ and $\sigma(X_t)$ are almost surely locally integrable ${\bf F}$-progressively measurable processes. 

The model (\ref{290920-1}) can express a fairly general class of models. %, as commented in Remark \ref{290920-2}. 
For example, consider a system
\beas 
Y_t &=& Y_0+\int_0^t b(\gamma,s,L_s,\xi_s,\theta^*)ds+\int_0^t\sigma(\gamma,s,L_s,\xi_s)dw_s\\
\xi_t &=& \xi_0+\int_0^t \tilde{b}(\xi_s)ds+\int_0^t \tilde{\sigma}(\xi_s)d\tilde{w}_s
\eeas
where $\gamma$ expresses a measurable random scenery taking values in a measurable space $({\sf G},\bbB_{\sf G})$, and 
$(b(\gamma,\cdot),\sigma(\gamma,\cdot))$ are regarded as a random environment in space-time. 
$\xi_t$ is a latent diffusion process having a good mixing property. 
The process $L_t$ is a process with long memory. % but independent of $(\xi_t)$. 
The process $(Y_t,\xi_t)$ is like a diffusion process but it does not enjoy a fast decay of mixing coefficient due to the component $L_t$. 
In this example, we may set $X_s=(\gamma,s,L_s,\xi_s)$. 
It is also possible to incorporate feedback of $Y_t$ as $X_s=(\gamma,s,L_s,\xi_s,Y_s)$. 
If the whole path $(L_t)$ is included in $\gamma$, then a simplified expression $b(\gamma,s,\xi_s,\theta^*)$ is possible 
for $b(\gamma,s,L_s,\xi_s,\theta^*)$.

We estimate the true value $\theta^*$ of the parameter $\theta\in\Theta$ based on observations $((Y_t,X_t)_{t\in[0,T]})$. 
Let $S=\sigma\sigma^\star$ and assume that $S(X_t)$ is invertible a.s. 
Define a random function $\bbH_T$ by 
\beas 
\bbH_T(\theta) &=& 
\int_0^T S(X_t)^{-1}\big[b(X_t,\theta),dY_t\big]
-\half \int_0^T S(X_t)^{-1}\big[b(X_t,\theta)^{\otimes2}\big]dt.
\eeas
By (\ref{290920-1}), $\bbH_T$ has the following representation: 
\beas 
\bbH_T(\theta) 
&=&
M_T(\theta)+N_T(\theta),
\eeas
where 
\beas 
M_T(\theta) &=& 
\int_0^T S(X_t)^{-1}\big[b(X_t,\theta),\sigma(X_t)dw_t\big]
\eeas
and 
\beas 
N_T(\theta) &=& 
\int_0^T S(X_t)^{-1}\bigg[b(X_t,\theta)\otimes b(X_t,\theta^*)-\half b(X_t,\theta)^{\otimes2}\bigg]dt.
\eeas

Define a $(\sfr+1)$-dimensional function $H$ by 
\beas 
H(x,\theta) &=& \bigg(S(x)^{-1}\big[b(x,\theta),\sigma(x)\cdot\big],\>
S(x)^{-1}\bigg[b(x,\theta)\otimes b(x,\theta^*)-\half b(x,\theta)^{\otimes2}\bigg]\bigg)
\eeas
for $x\in{\sf X}$ and $\theta\in\Theta$.

\bd
\im[[C1\!\!]] %{\bf (i)} 
%\bd
%\im[(ii)] 
The mapping $\Theta\ni\theta\mapsto H(x,\theta)$ is four times continuously differentiable and 
\beas 
\sup_{\theta\in\Theta}\sum_{i=0}^4\big|\partial_\theta^i H(x,\theta)\big| &\leq& H_1(x)
\quad(x\in{\sf X})
\eeas
for some measurable function 
$H_1:{\sf X}\to\bbR_+$ such that 
\beas 
\sup_{t\in\bbR+}E\big[H_1(X_t)^p\big]<\infty\quad a.s. 
\eeas
for all $p>1$. 
%\koko
%\im[(iii)] $\ds E_\calc\big[\partial_\theta H(\gamma,X_t,\theta^*)\big]=0$ a.s. for $t\in\bbR_+$. 
%\ed
\ed\halflineskip

Let $\calc$ be a sub $\sigma$-field of $\calf_0$, and let 
$\calb_I=\calc\vee\sigma[X_t,w_t-w_{\inf I};\> t\in I]$ for $I\subset\bbR_+$. 
A partial mixing coefficient $\alpha_\calc(h)$ is a $\calc$-measurable 
$[0,1/2]$-valued random variable satisfying the inequality
\beas 
\alpha_\calc(h) &\geq& 
\sup_{t\in\bbR_+}\sup\bigg\{\big|P_\calc\big[A\cap B\big]-P_\calc\big[A\big]P_\calc\big[B\big]\big|;\>
A\in\calb_{[0,t]},\>B\in\calb_{[t+h,\infty)}\bigg\}
\eeas
on some $\Omega_0\in\calf$ such that $P[\Omega_0]=1$. 
%Assume that $\calc\subset\calf_0$. 

{\colred Let ${\sf X}=\bbR^\sfd$.} 
Suppose that a regular conditional probability $\mu_t=P_\calc^{X_t}[\cdot]$ of $X_t$ given $\calc$ exists. 
The measure-valued process $\mu_t$ is a ``basso continuo'', which may only admit 
a very weak ergodic property. 
We will consider the following two situations. 

\bd
\im[[C2\!\!]] {\bf (i)} There exists a positive constant $L_0$ such that for every $L>0$, 
\beas 
\limsup_{h\to\infty} h^L\big\|\alpha_\calc(h)^{L_0}\big\|_1&<&\infty.
\eeas 
\bd\im[(ii)] 
There exist a probability measure $\nu$ on $\bbR^\sfd$ and a positive constant $\ep_1$ such that 
\beas 
T^{\ep_1}\bigg|\frac{1}{T}\int_0^T\mu_t(f)dt - \nu(f)\bigg|&\to&0\quad a.s.
\eeas
as  $T\to\infty$ for %any bounded measurable function $f:\bbR^\sfd\to\bbR$. 
any measurable function $f:\bbR^\sfd\to\bbR$ satisfying 
$|f(x)|\leq C(1+H_1(x)^C)$ ($x\in\bbR^\sfd$) for some positive constant $C$. 
\ed
\ed
\halflineskip

Here we wrote $\mu_t(f)=\int f(x)\mu_t(dx)$ for a measurable function $f:\bbR^\sfd\to\bbR$. 
The strong mixing coefficient of the measure valued process $\mu=(\mu_t)_{t\in\bbR_+}$ is defined by 
\beas 
\alpha^\mu(h) &=& \sup\bigg\{\big|P[A\cap B]-P[A]P[B]\big|;\>A\in \calc_{[0,t]},\> B\in\calc_{[t,\infty)}\bigg\}
\eeas
where 
$\calc_I=\sigma\big[\mu_t(f);\>t\in I,\>f\in C_b(\bbR^\sfd)\big]$ for $I\subset\bbR_+$. 

\bd
\im[[C2$^\sharp$\!\!]] {\bf (i)} There exists a positive constant $L_0$ such that for every $L>0$, 
\beas 
\limsup_{h\to\infty} h^L\big\|\alpha_\calc(h)^{L_0}\big\|_1&<&\infty.
\eeas 
\bd\im[(ii)] 
For some $\ep_0>0$, %$\sum_{h=1}^\infty h^{-1+\ep_0}\alpha^\mu(h)<\infty$. 
$\alpha^\mu(h)=O(h^{-\ep_0})$ as $h\to\infty$. 
\im[(iii)] There exist a probability measure $\nu$ on $\bbR^\sfd$ and a positive constant $\ep_1$ such that 
\beas 
T^{\ep_1}\bigg|\frac{1}{T}\int_0^TE\bigg[\int_{\bbR^\sfd}f(x)\mu_t(dx)\bigg]dt - \nu(f)\bigg|&\to&0
\eeas
as  $T\to\infty$ 
%for any bounded measurable function $f:\bbR^\sfd\to\bbR$. 
for %any bounded measurable function $f:\bbR^\sfd\to\bbR$. 
{\colorr 
any measurable function $f:\bbR^\sfd\to\bbR$ satisfying 
$|f(x)|\leq C(1+H_1(x)^C)$ ($x\in\bbR^\sfd$) for some constant $C$.} 
\ed
\ed
\halflineskip

\begin{en-text}
Let 
\beas 
\bbY(\theta) &=& E\big[H(L_0,X_0,\theta)-H(L_0,X_0,\theta^*)\big]. 
\eeas
\end{en-text}

Define $\bbY:\Theta\to\bbR$ by 
\beas 
\bbY(\theta) 
&=&
-\half \int_{\bbR^\sfd}
S(x)^{-1}\big[\big(b(x,\theta)-b(x,\theta^*)\big)^{\otimes2}\big]\nu(dx).
\eeas

\bd
\im[[C3\!\!]] %{\bf (i)} 
%$\bbY(\theta)<0$ for any $\theta\in\bar{\Theta}\setminus\{\theta^*\}$. 
There exists a positive constant $\chi_0$ such that 
$\bbY(\theta)\leq -\chi_0|\theta-\theta^*|^2$ for all $\theta\in\Theta$. 
\ed
\halflineskip

Under $[C3]$, the matrix 
\beas 
\Gamma &:=& -\partial_\theta^2\bbY(\theta^*)
\yeq 
\int_{\bbR^\sfd}
S(x)^{-1}\big[\big(\partial_\theta b(x,\theta^*)\big)^{\otimes2}\big]\nu(dx)
\eeas
is a positive-definite $\sfp\times\sfp$ symmetric matrix. 
%The equality $\partial_\theta\bbY(\theta^*)=E\big[\partial_\theta H(L_0,X_0,\theta^*)\big]=0$ also follows but 
%we have assumed a stronger property $E_\calc\big[\partial_\theta H(L_t,X_t,\theta^*)\big]=0$. 
Let $\hat{\theta}_n^M=\hat{\theta}_T$ and $\hat{\theta}_n^B=\tilde{\theta}_T$, and let 
$\hat{u}_n^M=\hat{u}_T$ and $\hat{u}_n^B=\tilde{u}_T$. 

\halflineskip
\begin{theorem}\label{290922-2}
{\bf (i)} 
Suppose that Conditions $[C1]$, $[C2]$ and $[C3]$ are satisfied. Then 
\bea\label{291009-1}
\hat{u}^{\sf A}_T\to^d\Gamma^{-1/2}\zeta
\eea
as $T\to\infty$ for ${\sf A}=M$ and $B$, 
where $\zeta$ is a $\sfp$-dimensional standard Gaussian random vector. 
\bd\im[(ii)] 
The convergence (\ref{291009-1}) holds under Conditions $[C1]$, $[C2^\sharp]$ and $[C3]$. 
\ed
\end{theorem}
\halflineskip

%\noindent
\subsection{Proof of Theorem \ref{290922-2}} 
Let $\ep_*>0$. 
Define $\Psi_T$ by 
\beas 
\Psi_T &=& 1_{A_T}
\eeas
where 
\beas 
A_T &=& \left\{
\max_{j=1,...,\lceil T\rceil}E_\calc\bigg[\int_{j-1}^jH_1(X_t)^{r_*}dt\bigg]\leq T^{\ep_*}
\right\},
\eeas
where we fix a sufficiently large but finite constant $r_*$ 
in what follows, since we only aim at asymptotic normality of the QLA estimators. 

Let 
$\bbY_T(\theta) = T^{-1}\big(\bbH_T(\theta)-\bbH_T(\theta^*)\big)$. 
Then
\beas 
\bbY_T(\theta) 
&=& 
\bbY_T^{(0)}(\theta)+\bbY_T^{(1)}(\theta),
\eeas
where 
\beas 
\bbY_T^{(0)}(\theta) &=& 
\frac{1}{T}\big\{M_T(\theta)-M_T(\theta^*)\big\}
\eeas
and 
\beas 
\bbY^{(1)}_T(\theta)&=&
-\frac{1}{2T}\int_0^TS(X_t)^{-1}\bigg[\big(b(X_t,\theta)-b(X_t,\theta^*)\big)^{\otimes2}\bigg]dt.
\eeas

Let 
\beas 
\Delta_T&=&T^{-1/2}\partial_\theta\bbH_T(\theta^*)
\yeq T^{-1/2}\partial_\theta M_T(\theta^*) 
\\&=& 
T^{-1/2}\int_0^T S(X_t)^{-1}\big[\partial_\theta b(X_t,\theta^*),\sigma(X_t)dw_t\big]
\eeas
and let
\beas 
\Gamma_T&=&-T^{-1}\partial_\theta^2\bbH_T(\theta^*)
\yeq 
-T^{-1}\partial_\theta^2 M_T(\theta^*) -T^{-1}\partial_\theta^2N_T(\theta^*) 
\\&=&
-T^{-1}\int_0^T S(X_t)^{-1}\big[\partial_\theta^2 b(X_t,\theta^*),\sigma(X_t)dw_t\big]
\\&&
+T^{-1}\int_0^T S(X_t)^{-1}\big[\big(\partial_\theta b(X_t,\theta^*)\big)^{\otimes2}\big]dt.
\eeas

\begin{lemma}\label{290919-1}
{\bf (i)} $\ds \sup_{T\in\calt}E_\calc\big[|\Delta_T|^M\big]%{\colorr \Psi_T}
<\infty$ a.s. for every $\calt\in{\mathfrak T}$ and every $M>0$. 
%{\colorr $L_t$の（最大値ではなく）平均量で抑えられないと，ここが通らない．$\bar{H}$ bounded（あるいは$L_t$部分の影響が有界）しか通らない．}
\bd\im[(ii)] Let $\eta\in(0,1/2)$ and $M>0$. Then for sufficiently large $r_*$, one has
\beas 
\sup_{T\in\calt}E_\calc\bigg[\bigg(\sup_{\theta\in\Theta} T^{\half-\eta}\big|\bbY_T(\theta)-\bbY(\theta)\big|\bigg)^M\bigg]{\colorr \Psi_T}<\infty\quad a.s. 
\eeas
for every $\calt\in{\mathfrak T}$. %and every $M>0$. 
\ed
\end{lemma}
\proof 
%Write $\mu_t(g)=\int g(x)\mu_t(dx)$ for a measurable function $g:\bbR^\sfd\to\bbR$. 
Suppose that $|g(x)|\leq C(1+H_1(x)^C)$ for some constant $C>0$. 
Then 
\bea\label{290921-2} 
\sup_{t\in\bbR_+}E\big[|\mu_t(g)|^p\big] 
&\leq&
\sup_{t\in\bbR_+}E\big[|g(X_t)|^p\big]
%\>\leq\>
%\sup_{t\in\bbR_+}E\big[H_1(X_t)^{p}\big]
\><\>\infty
\eea
for any $p\geq1$. 
\begin{en-text}
Let $\ep>0$. 
Then 
\beas 
\sup_{t\in\bbR_+}E\bigg[\int|g(x)|1_{\{|g(x)|\geq C\}}\mu_t(dx)\bigg] 
&\leq&
\sup_{t\in\bbR_+}\frac{1}{C}E\bigg[\int|g(x)|^2\mu_t(dx)\bigg] 
\\&\leq&
\sup_{t\in\bbR_+}\frac{1}{C}E\big[H_1(X_t)^{2p}]
\><\>
\ep
\eeas
for sufficiently large $C$. 
\end{en-text}

{\colorr Suppose that $[C2^\sharp]$ holds, for a while.}
Let $r\in(1,\min\{1+\ep_0,2\})$. 
In the notation of Rio \cite{rio2017asymptotic}, for the tail-quantile function $Q_j(u)$ of $\int_{j-1}^j\mu_t(g)dt$, we have 
$Q_j(u)\simleq u^{-1/L}$ for arbitrarily large $L$ due to $L_{\infty-}$-boundedness of $H_1(X_t)$ uniform in $t$. 
Moreover, $(\alpha^\mu)^{-1}(u)\simleq u^{-1/\ep_0}$. Therefore, 
\beas 
M_{r,\alpha^\mu}(Q)
&=& 
\int_0^1 \big[(\alpha^\mu)^{-1}(u)\big]^{r-1}Q(u)^rdu 
\>\simleq\> 
\int_0^1 u^{^{-\ep_0^{-1}(r-1)-L^{-1}r}}du
\><\> \infty
\eeas
if we take a sufficiently large $L$. 
We apply Corollary 3.2 (i) of Rio \cite{rio2017asymptotic} to conclude 
\beas
\frac{1}{n^{1/r}}\bigg(\int_0^n\mu_t(g)dt-\int_0^nE\big[\mu_t(g))\big]dt \bigg)
&\to&
0\quad a.s.
\eeas
as $T\to\infty$. 
Further, since 
\beas
\sum_{n\in\bbN}P\bigg[\sup_{T:n\leq T<n+1}\int_n^{T}|\mu_t(g)|dt>n^{1/(2r)}\bigg]&<&\infty, 
\eeas
we obtain
\bea\label{290921-1} 
\frac{1}{T^{1/r}}\bigg(\int_0^T\mu_t(g)dt-\int_0^TE\big[\mu_t(g)\big]dt \bigg)
&\to&
0\quad a.s.
\eea
as $T\to\infty$. 
\begin{en-text}
Let $\ep_1\in(0,\ep_0)$. 
\beas&&
\bigg|\frac{1}{n}\int_0^n\int g(x)\mu_t(dx)dt-\frac{1}{n}\int_0^nE\bigg[\int g(x)\mu_t(dx)\bigg]dt \bigg|
\\&=&
\bigg|\frac{1}{n}\sum_{j=1}^{n}\int_{j-1}^{n}\bigg(\int g(x)\mu_t(dx)-E\bigg[\int g(x)\mu_t(dx)\bigg]\bigg)dt \bigg|
\\&\leq&
\bigg|\frac{1}{n^{1-\ep_1}}\sum_{j=1}^{n}\int_{j-1}^{j}\bigg(\int n^{-\ep_1}g(x)1_{\{|g(x)|\leq n^{\ep_1}\}}\mu_t(dx)
-E\bigg[\int n^{-\ep_1}g(x)1_{\{|g(x)|\leq n^{\ep_1}\}}\mu_t(dx)\bigg]\bigg)dt \bigg|
\\&&
+\frac{1}{n}\sum_{j=1}^{n}\int_{j-1}^{j}\bigg(\int |g(x)|1_{\{|g(x)|> n^{\ep_1}\}}\mu_t(dx)
+E\bigg[\int |g(x)|1_{\{|g(x)|>n^{\ep_1}\}}\mu_t(dx)\bigg]\bigg)dt 
\eeas
\end{en-text}
We choose a sufficiently large constant $\eta\in(0,1/2)$. 
Under $[C2^\sharp]$ (iii) %by cut-off argument 
with (\ref{290921-2}), we have 
\beas 
T^{\half-\eta}\bigg|\frac{1}{T}\int_0^T E\big[\mu_t(g)\big]dt - \nu(g)\bigg| &\to& 0
\eeas
as $T\to\infty$. Then (\ref{290921-1}) gives 
\bea\label{290921-3}
T^{\half-\eta}\bigg|\frac{1}{T}\int_0^T \mu_t(g)dt - \nu(g)\bigg|&\to&0\quad a.s.
\eea
as $T\to\infty$. 
Under $[C2]$, the convergence (\ref{290921-3}) is obvious for a suitable $\eta$. 

By the Burkholder-Davis-Gundy inequality, 
\bea\label{290921-4} 
E_\calc\big[|\Delta_T|^M\big]
&\simleq& 
E_\calc\bigg[\bigg|\frac{1}{T}\int_0^TH_1(X_t)^2dt\bigg|^{M/2}\bigg]
\nn\\&\simleq&
\frac{1}{T}\int_0^T\mu_t(H_1^{M})dt
\>\to\>
\nu(H_1^{M})\quad a.s.
\eea
as $T\to\infty$ for $M\geq2$, which proved (i). 

In this situation, we can exterchange the differentiation in $\theta$ and the stochastic integral, and 
\beas 
\partial_\theta\bbY^{(0)}_T(\theta) 
&=& 
\frac{1}{T}\int_0^T S(X_t)^{-1}\big[\partial_\theta b(X_t,\theta),\sigma(X_t)dw_t\big].
\eeas
Then with Sobolev's inequality and following the way in (\ref{290921-4}), we obtain 
\beas 
E_\calc\big[\|T^{2^{-1}-\eta}\>\bbY_T^{(0)}\|_{C(\Theta)}^M\big]
&\simleq&
\int_\Theta \big\{E_\calc\big[|T^{2^{-1}-\eta}\>\bbY_T^{(0)}(\theta)|^M]+E_\calc\big[|T^{2^{-1}-\eta}\>\partial_\theta\bbY_T^{(0)}(\theta)|^M]\big\}d\theta
\\&\to&
0\quad a.s.
\eeas
as $T\to\infty$. 

Next, we apply Lemma \ref{290824-1} to $\bbY_T^{(1)}(\theta)-\bbY(\theta)$ and 
$\partial_\theta\big(\bbY_T^{(1)}(\theta)-\bbY(\theta)\big)$ 
with the help of $\Psi_T$, as well as Sobolev's inequality, 
to show (ii). More precisely, let 
\beas 
g(x,\theta) &=& -\half S(x)^{-1}\big[(b(x,\theta)-b(x,\theta^*))^{\otimes2}\big]. 
\eeas
Then {\colred for $M>\sfp$,}
\beas&&
E_\calc\big[\|T^{2^{-1}-\eta}\>\big(\bbY_T^{(1)}-\bbY\big)\Psi_T\|_{C(\Theta)}^M\big]
\\&\simleq&
\sum_{i=0,1}
\int_\Theta E_\calc\big[|T^{2^{-1}-\eta}\>\partial_\theta^i\big(\bbY_T^{(1)}(\theta)-\bbY(\theta)\big)\Psi_T|^M]d\theta
\\&\simleq&
I_T+J_T
\eeas
where 
\beas 
I_T &=& 
\sum_{i=0,1}
\int_\Theta E_\calc\big[|T^{2^{-1}-\eta}\>\big(\partial_\theta^i\bbY_T^{(1)}(\theta)
-E_\calc\big[\partial_\theta^i\bbY_T^{(1)}(\theta)\big]\big)\Psi_T|^M]d\theta
\eeas
and 
\beas 
J_T &=& 
\sum_{i=0,1}
\int_\Theta |T^{2^{-1}-\eta}\>\big(E_\calc\big[\partial_\theta^i\bbY^{(1)}_T(\theta)\big]-\partial_\theta^i\bbY(\theta)\big)\Psi_T|^Md\theta. 
%\\&\to&
%0\quad a.s.
\eeas
We notice that $\nu(H_1^p)<\infty$ for every $p>1$ by (\ref{290921-2}) 
under both $[C2]$ and $[C2^\sharp]$; in particular, 
$\partial_\theta\bbY(\theta)=\nu(\partial_\theta g(\cdot,\theta))$. 
By (\ref{290921-3}), we have 
$J_T\to0$ as $T\to\infty$ a.s. 
Applying Lemma \ref{290824-1}, we see 
\beas 
I_T &\leq&
T^{M(-\eta+\frac{\ep_*}{r_*})}V_*
\eeas
for suitably set $(M,r_*)$ so that $r_*>M\geq2$ and $-\eta+\ep_*/r_*<0$, where 
\beas 
V_* &\leq& C(M,r_*)
\bigg[\bigg(1+\sum_{h=1}^{\infty}\alpha_\calc(h)^{1-\frac{2}{r_*}}\bigg)^{\colorr M/2}+\sum_{h=1}^\infty(h+1)^{M-2}\alpha_\calc(h)^{1-\frac{M}{r_*}}\bigg]. 
\eeas
\begin{en-text}
We can take sufficiently large constants $M$, $r_*$ and $K$ so that 
\beas 
r_*>M\geq2,\>\left(1-\frac{M}{r_*}\right)K\geq L_0;
\eeas
in particular,  $r_*>\frac{M}{M-1}$. 
\end{en-text}
Let $K=L_0(1-M/r_*)^{-1}$. 
Since $\alpha_\calc(h)\leq1/2\leq1$, we have 
\beas 
\bigg\|\bigg(\sum_{h=1}^{\infty}\alpha_\calc(h)^{1-\frac{2}{r_*}}\bigg)^K\bigg\|_1
&\leq&
\bigg\|\bigg(\sum_{h=1}^{\infty}(h+1)^{M-2}\alpha_\calc(h)^{1-\frac{M}{r_*}}\bigg)^K\bigg\|_1
\\&\leq&
\bigg\|\sum_{h=1}^\infty (h+1)^{-2}\bigg((h+1)^M\alpha_\calc(h)^{1-\frac{M}{r_*}}\bigg)^K\bigg\|_1
\\&\leq&
\sum_{h=1}^\infty (h+1)^{KM-2}\big\|\alpha_\calc(h)^{L_0}\big\|_1
\><\>\infty. 
\eeas
Therefore $V_*<\infty$ a.s. and hence $I_T\to0$ as $T\to\infty$ a.s.  
\qed\halflineskip

In a similar fashion to Lemma \ref{290919-1}, we can show the following lemma. 

\begin{lemma}\label{290919-2} 
{\bf (i)} 
Let $M>0$. Then for a sufficiently large $r_*$, for any $\calt\in{\mathfrak T}$, %and $M>0$, 
\beas 
\sup_{T\in \calt}E_\calc\bigg[\bigg(T^{-1}\sup_{\theta\in\Theta}
\big|\partial_\theta^3\bbH_T(\theta)\big|\bigg)^{M}\bigg]{\colorr \Psi_T}
&<&\infty\quad a.s.
\eeas
\bd\im[(ii)] 
Let $M>0$ and 
$\eta\in(0,1/2)$. Then for a sufficiently large $r_*$, 
\beas 
\sup_{T\in \calt}E_\calc\bigg[\big(T^{\eta}\big|\Gamma_T-\Gamma\big|\big)^{M}\bigg]{\colorr \Psi_T}
&<&\infty\quad a.s.
\eeas
for any $\calt\in{\mathfrak T}$.
\ed
\end{lemma}
\halflineskip

As before, the random field $\bbZ_T$ is defined by 
\beas 
\bbZ_T(u) 
&=& 
\exp\big(\bbH_T(\theta^\dagger_T(u))-\bbH_T(\theta^*)\big)\quad (u\in\bbR^\sfp)
\eeas

\begin{lemma}\label{290919-3}
%Conditional PLD. 
Let $L>0$. Then there exist $r_*>0$ (in $\Psi_T$) and $\varrho\in(1,2)$ such that 
\beas
\sup_{T\in\calt}\sup_{r\in\bbN}r^LP_\calc\bigg[\sup_{u\in\bbV_T(r)}\bbZ_T(u)\geq \exp\big(-2^{-1}r^{\varrho}\big)\bigg]\Psi_T
&<&
\infty\quad a.s.
\eeas
for every $\calt\in{\mathfrak T}$. 
\end{lemma}
\proof 
Apply Theorem \ref{290916-1} with the help of Lemmas \ref{290919-1} and \ref{290919-2}. 
\qed\halflineskip

\begin{lemma}\label{290919-4}
%Modulus of continuity
For any $\ep>0$ and $c>0$, 
$\ds \limsup_{T\to\infty}P_\calc\big[W_T(\delta,c,\ep)\big]\Psi_T\to^P0$ as $\delta\down0$. 
\end{lemma}
\proof
We apply Lemma \ref{290919-2} (i) to estimate $r_T(u)$ in the LAQ representation of $\bbZ_T$. 
\begin{en-text}
In order to estimate $r_T(u)$, 
it suffices to show that %$\limsup_{T\to\infty}E\big[1_{W_T(\delta,c,\ep)}\Psi_T\big]\to0$ as $\delta\down0$. 
\beas 
\limsup_{T\to\infty,T\in\calt}\sum_{i=2}^4\sup_{\theta\in\Theta}E_\calc\big[\big|\partial_\theta^i \bbY_T(\theta)\big|^{\sfp+1}\big]\Psi_T
&<&\infty
\quad a.s.
\eeas
for every $\calt\in{\mathfrak T}$. 
The Burkholder-Davis-Gundy inequality and Lemma \ref{290824-1} serve it;  
as a matter of fact, the procedures are essentially the same as those in the proof of 
Lemmas \ref{290919-1} and \ref{290919-2}. 
\end{en-text}
\qed\halflineskip

Define a random field $\bbZ:\overline{\Omega}\times\bbR^\sfp\to\bbR$ on an extension $(\overline{\Omega},\overline{\calf},\overline{P})$ 
of $(\Omega,\calf,P)$ by 
\beas 
\bbZ(u) 
&=& 
\exp\bigg(\Delta[u]-\half\Gamma[u^{\otimes2}]\bigg),
\eeas
where 
$\Delta=\Gamma^{1/2}\zeta$ and $\zeta$ is a $\sfp$-dimensional standard Gaussian random variable defined on $\overline{\Omega}$ 
and independent of $\calf$. 

\begin{lemma}\label{290919-5}
%Conditional CLT
{\bf (i)} 
For any $k\in\bbN$, $u_i\in\bbR^\sfp$ $(i=1,...,k)$ and $f\in C_b(\bbR^{k\sfp})$, 
\beas 
E_\calc\big[f\big((\bbZ_T(u_i))_{i=1,...,k}\big)\Psi_T\big]
&\to^P& 
\overline{E}_\calc\big[f\big((\bbZ(u_i))_{i=1,...,k}\big)\big]\yeq \overline{E}\big[f\big((\bbZ(u_i))_{i=1,...,k}\big)\big]
\eeas
as $T\to\infty$. 
\bd
\im[(ii)] 
$\Psi_T\to^P1$ as $T\to\infty$. 
\ed
\end{lemma}
\proof 
The conditional version of martingale central limit theorem gives 
\bea\label{291010-3}
E_\calc[g(\Delta_T)]\to^P \overline{E}_\calc[g(\Delta)]
\eea
for $g\in C_b(\bbR^\sfp;\bbR^k)$.
Indeed, the quadratic variation of the martingale associated with $\Delta_T$ is 
$\frac{1}{T}\int_0^T g(X_t,\theta^*)^{\otimes2}dt$ if evaluated at $T$, and 
it converges to $\Gamma$ in probability. 
Then we have the convergence $E_\calc\big[\Psi_T\exp\big(\Delta_T[iu]+2^{-1}\Gamma[u^{\otimes2}]\big)\big]\to^p1$ as $T\to\infty$ for every $u\in\bbR^\sfp$. 
We obtain (\ref{291010-3}) with uniform approximation of $g$ on a compact set  
by trigonometric functions. 

\begin{en-text}
\koko by stability of convergence, one has 
\beas 
\big\|E_\calc[g(\Delta_T)]-\overline{E}_\calc[g(\Delta)]\big\|_{L^2(\calc)}
&=&
\sup\bigg\{E\bigg[\big(E_\calc[g(\Delta_T)]-\overline{E}_\calc[g(\Delta)]\big)C\bigg];\> 
C\in L^2(\calc)\bigg\}
\\&=&
\sup\bigg\{E[g(\Delta_T)C]-\overline{E}[g(\Delta)C];\> 
C\in L^2(\calc)\bigg\}
\eeas
 for every bounded $\calc$-measurable function $C$, 
that is, $E_\calc[g(\Delta_T)]-\overline{E}_\calc[g(\Delta)]\to 0$ weakly in $L^1(\overline{\Omega},\overline{\calf},\overline{P})$. 
\end{en-text}
In the representation (\ref{290922-1}) of $\bbZ_T$, the convergence 
$E_\calc[|r_T(u)|\wedge 1]\to^P0$ follows from e.g. Lemma \ref{290919-2} (i). 
Thus we obtain (i). The property (ii) is easy to show by definition of $\Psi_T$. 
\qed\halflineskip

Condition $[A5]$ is verified e.g. with Lemma 2 of \cite{yoshida2011polynomial}. 
Now Theorem \ref{290922-2} follows from Theorems \ref{170902-2} and \ref{170904-3} 
together with 
{\colred Theorem \ref{290916-1} as well as} 
Lemmas \ref{290919-3}, \ref{290919-4} and \ref{290919-5}. 
%\qed\halflineskip

\subsection{An example}
On a stochastic basis $(\Omega,\calf,{\bf F},P)$, ${\bf F}=(\calf_t)_{t\in\bbR_+}$, let us consider stochastic processes 
$Y=(Y_t)_{t\in\bbR_+}$, $L=(L_t)_{t\in\bbR_+}$ and $U=(U_t)_{t\in\bbR_+}$ satisfying 
\beas 
Y_t &=& Y_0+\int_0^t b_0(L_s)b_1(U_s,\theta^*)ds+\int_0^t\sigma_0(L_s)\sigma_1(U_s)dw_s
%L_t &=& L_0-a\int_0^t L_sds+dB^H_s. \koko
%U_t &=& U_0+\int_0^t \tilde{b}(U_s)ds+\int_0^t \tilde{\sigma}(U_s)d\tilde{B}_s
\eeas
where $w=(w_t)_{t\in\bbR_+}$ is a one-dimensional standard ${\bf F}$-Wiener process. 
%Here $B=(B_t)_{t\in\bbR_+}$ is a Brownian motion and $B^H=(B^H_t)_{t\in\bbR_+}$ is a fractional Brownian motion with Hurst parameter $H$. 
We assume 
\bd\im[(i)] %$a$ is a positive constant. 
$L$ is \cadlag $\calf_0$-measurable, stationary and independent of $(U,w,Y_0)$. 
The $\alpha$-mixing coefficient $\alpha^L$ of $L$ satisfies $\alpha^L(h)\simleq h^{-a}$ as $h\to\infty$ for some positive constant $a$. 
For every  $p>1$, %$\sup_{t\in\bbR_+}\|L_t\|_p<\infty$. 
$\|L_0\|_p<\infty$. 
\im[(ii)] $U$ is a \cadlag ${\bf F}$-progressively measurable stationary process satisfying 
%$\sup_{t\in\bbR_+}\|U_t\|_p<\infty$ for every $p>1$. 
$\|U_0\|_p<\infty$ for every $p>1$. 
The $\alpha$-mixing coefficient $\alpha^U$ of $U$ satisfies $\alpha^U(h)\leq b^{-1}e^{-bh}$ for some positive constant $b$. 
\im[(iii)] %$b_0$ and $\sigma_0$ are measurable functions such that $\lim_{|x|\to\infty}\log(1+|b_0(x)|+|\sigma_0(x))/|x|^2=0$. 
$b_0$, $\sigma_0$ and $\sigma_1$ are measurable functions of at most polynomial growth. 
The function 
$b_1$ is measurable, four times continuously differentiable in $\theta$ and $\partial_\theta^ib_1(\cdot,\theta)$ is of at most polynomial growth 
uniformly in $\theta\in\Theta$ for $i=0,...,4$. 
Moreover $\inf_{z,u}\sigma_0(z)\sigma_1(u)>0$. 
\ed

The variable $X_t=(L_t,U_t)$ for this model. 
The random field $\bbH_T$ is given by 
\beas 
\bbH_T(\theta)
&=& 
\int_0^T \frac{b_0(L_t)b_1(U_t,\theta)}{\sigma_0(L_t)^2\sigma_1(U_t)^2}dY_t
-\half\int_0^T \frac{b_0(L_t)^2b_1(U_t,\theta)^2}{\sigma_0(L_t)^2\sigma_1(U_t)^2}dt. 
\eeas
It has a representation 
\beas 
\bbH_T(\theta)
&=& 
\int_0^T \frac{b_0(L_t)b_1(U_t,\theta)}{\sigma_0(L_t)\sigma_1(U_t)}dw_t
+\int_0^T 
\bigg\{\frac{b_0(L_t)^2b_1(U_t,\theta)b_1(U_t,\theta^*)}{\sigma_0(L_t)^2\sigma_1(U_t)^2}
-\half\frac{b_0(L_t)^2b_1(U_t,\theta)^2}{\sigma_0(L_t)^2\sigma_1(U_t)^2}\bigg\}dt. 
\eeas
Let $\calc=\sigma[L_t;\>t\in\bbR_+]$. 
Since $\calb_I\equiv\calc\vee\sigma[L_t,U_t,w_t-w_{\inf I};\>t\in I]=\calc\vee\sigma[U_t,w_t-w_{\inf I};\>t\in I]$ 
for $I\subset\bbR_+$ and $\calc$ is independent of $\sigma[U_t,w_t;\>t\in\bbR_+]$, we can take $\alpha_\calc(h)=\alpha^{U,dw}(h)$, 
which is the $\alpha$-mixing coefficient associated with $\calb^{U,dw}_I=\sigma[U_t,w_t-w_{\inf I};\>t\in I]$ for $I\subset \bbR_+$. 
The coefficient $\alpha^{U,dw}$ enjoys an exponential decay; see Kusuoka and Yoshida \cite{KusuokaYoshida2000}. 
%\beas 
%P_\calc\big[g(L,U,w)f(L,U,w)\big]-P_\calc\big[g(L,U,w)\big]P_\calc\big[f(L,U,w)\big]
%\eeas

For any bounded measurable function $f$ on $\bbR^2$, 
\beas 
\mu_t(f)&=&
E_\calc\big[f(L_t,U_t)\big]
\yeq 
E\big[f(\ell,U_t)\big]\big|_{\ell=L_t}\yeq 
E\big[f(\ell,U_0)\big]\big|_{\ell=L_t}\quad a.s.
\eeas
In particular, $[C2^\sharp]$ (iii) holds for 
$\nu(f) = E\big[f(L_0,U_0)\big]$. 
Moreover, 
\beas 
\bbY(\theta) &=& 
-\half E\bigg[\frac{b_0(L_0)^2\big(b_1(U_0,\theta)-b_1(U_0,\theta^*)\big)^2}
{\sigma_0(L_0)^2\sigma_1(U_0)^2}\bigg]
\eeas
Thus, if $[C3]$ is satisfied, then $\hat{u}^{\sf A}_T$ (${\sf A}=M,B$) are asymptotically normal with 
variance 
\beas 
\Gamma &=& E\bigg[\frac{b_0(L_0)^2\big(\partial_\theta b_1(U_0,\theta^*)\big)^2}
{\sigma_0(L_0)^2\sigma_1(U_0)^2}\bigg]. 
\eeas
\halflineskip

\section{Stochastic regression model for volatility in random environment}\label{290924-6}
Let $(\Omega',\calf',P')$ be a probability space and let $(\Omega'',\calf'',{\bf F})$ be a measurable space having a right-continuous filtration ${\bf F}=(\calf_t)_{t\in[0,T]}$.  
We consider a transition kernel $Q_{\omega'}(d\omega'')$ from $\Omega'$ to 
$(\Omega'',\calf'')$ 
The extension $(\Omega,\calf,P)$ of $(\Omega',\calf',P')$ is defined by 
$\Omega=\Omega'\times\Omega''$, $\calf=\calf'\times\calf''$ and 
$P(d\omega',d\omega'')=P'(d\omega')Q_{\omega'}(d\omega'')$. 
Let $\bbT=[0,T]$. We consider measurable processes 
$b:\Omega\times\bbT\to\bbR^\sfm$, 
$X:\Omega\times\bbT\to\bbR^\sfd$, 
$Y:\Omega\times\bbT\to\bbR^\sfm$ and 
$w:\Omega\times\bbT\to\bbR^\sfr$. 
A random variable $\gamma$ takes values in a measurable space 
$({\sf G},\bbB_{\sf G})$ defined on a probability space $(\Omega',\calf',P')$. 

For each $\omega'\in\Omega'$, on the stochastic basis 
$\calb_{\omega'}=(\Omega'',\calf'',{\bf F},Q_{\omega'})$, 
we suppose that 
$b(\omega',\cdot)=(b_t(\omega',\cdot))_{t\in[0,T]}$ and 
$X(\omega',\cdot)=(X_t(\omega',\cdot))_{t\in[0,T]}$ become progressively measurable processes, 
$w(\omega',\cdot)=(w_t(\omega',\cdot))_{t\in[0,T]}$ is an $\sfr$-dimensional ${\bf F}$-Wiener process, 
and 
the process $Y(\omega',\cdot)=(Y_t(\omega',\cdot))_{t\in[0,T]}$ satisfies 
%a stochastic regression model %given by a  stochastic integral equation
\bea\label{291011-1} 
Y_t &=& Y_0+\int_0^t b_sds+\int_0^t \sigma(\gamma,X_s,\theta^*)dw_s,\quad t\in[0,T]. 
\eea
Here the function $\sigma$ is an $\bbR^\sfm\otimes\bbR^\sfr$-valued %a $M(\sfd,\sfr;\bbR)$-valued 
measurable function defined on ${\sf G}\times\bbR^\sfd\times\Theta$ and $\Theta$ is a bounded domain in $\bbR^\sfp$, $\theta^*\in\Theta$, 
and 
we suppose 
the stochastic integrals on $\calb_{\omega'}$ in (\ref{291011-1}) are well defined. 
Detailed conditions for it will be specified below. 
%$b_s(\omega',\omega'')$ and $\sigma(\gamma(\omega'),X_s(\omega',\omega''),\theta)$ 
%are bounded on every compact set of $s$ for every $(\omega',\omega'')$ and $\theta\in\Theta$. 

We observe $(\gamma,(X_{t_j}, Y_{t_j})_{j=0,...,n})$, where $t_j=t^n_j=jT/n$, and want to estimate 
$\theta^*$ from the data. 
It is regarded that $\omega'$ denotes the state of a random environment. 
The variable $\gamma$ describes the partially observed state of the random environment. 
The process $b$ is assumed to be completely unobservable. 

The random field $\bbH_n$ is defined by 
\beas 
\bbH_n(\theta) 
&=& 
-\half\sum_{j=1}^n\bigg\{\log\det S(\gamma,X_{t_{j-1}},\theta)+h^{-1}S(\gamma,X_{t_{j-1}},\theta)^{-1}\big[\big(\Delta_jY\big)^{\otimes2}\big]\bigg\},
\eeas
where $h=T/n$ and $\Delta_jY=Y_{t_j}-Y_{t_{j-1}}$. 
The QMLE and QBE are defined by $\bbH_n$ as in Section \ref{290918-1}. 
The existence of continuous extension of $\bbH_n$ to $\bar{\Theta}$ for the QMLE, and 
the conditions for the prior density $\varpi$ of the QBE are assumed, as before. 
Moreover, we assume that the boundary of $\Theta$ is good as in Section \ref{290924-5}. 
%as good as it admits the ordinary Sobolev's inequality for 
%the embedding $W^{1,p}(\Theta)\inclusion C_b(\Theta)$ 
%for $p>\sfp$. 

Let $\calc=\calf'$ and denote by $E_\calc[V](\omega')$ the integral 
$\int_{\Omega''}V(\omega',\omega'')Q_{\omega'}(d\omega'')$ of 
a measurable function $V$ on $(\Omega,\calf)$. 
%Let $\calc=\sigma[\gamma]$, the $\sigma$-field generated by $\gamma$. 
Let $S=\sigma\sigma^\star$. 
We will work with the following conditions. 
Suppose that a $\calc$-measurable random variable 
$K_p:\Omega'\to\bbR_+$ is given for every $p\geq1$.  
\bd
\im[[D1\!\!]] {\bf (i)} 
\begin{en-text}
$E_\calc[|X_0|^p]\leq K_p$ %<\infty$ a.s. 
for all $p>0$. For every $p>0$, 
%there exists a $\calc$-measurable random varialbe $K_p$ such that 
\beas 
E_\calc\big[|X_t-X_s|^p\big]&\leq& K_p|t-s|^{p/2}
\eeas
for all $t,s\in[0,T]$. 
\end{en-text}
%%
%\bd
%\im[(ii)] 
$\sup_{t\in[0,T]}E_\calc\big[|b_t|^p\big]\leq K_p$. %<\infty$ a.s. for all $p>0$. 
\bd
\im[(ii)] 
The mapping $(x,\theta)\mapsto\sigma(\gamma,x,\theta)$ is continuously differentiable 
twice in $x$ and four times in $\theta$ and 
%there exists a nonnegative random variable $K$ such that 
\beas 
%\sup_{\gamma\in{\sf G}}
\sum_{i=0,1,2\atop j=0,1,2,3,4}|\partial_x^i\partial_\theta^j\sigma(\gamma,x,\theta)| &\leq& K_1(1+|x|)^{K_1}\quad 
\eeas
Furthermore, $\big(\inf_{x,\theta}\det S(\gamma,x,\theta)\big)^{-1}\leq K_1$. 
\im[(iii)] On each $\calb_{\omega'}$, the process $X(\omega',\cdot)=(X_t(\omega',\cdot))_{t\in[0,T]}$ admits a representation 
\beas 
X_t &=& X_0+\int_0^t \tilde{b}_sds+\int_0^t a_sdw_s+\int_0^t \tilde{a}_sd\tilde{w}_s,
\eeas
where $\tilde{w}=(\tilde{w}_t)_{t\in[0,T]}$ is an $\sfr_1$-dimensional ${\bf F}$-Wiener process independent of $w$, and 
$\tilde{b}=(\tilde{b}_t)_{t\in[0,T]}$, $a=(a_t)_{t\in[0,T]}$ and $\tilde{a}=(\tilde{a}_t)_{t\in[0,T]}$ are 
progressively measurable processes taking values in $\bbR^\sfd$, $\bbR^\sfd\otimes\bbR^\sfr$ and  $\bbR^\sfd\otimes\bbR^{\sfr_1}$, 
respectively, satisfying 
\beas
E_\calc\big[|X_0|^p\big]+\sup_{t\in[0,T]}\big(E_\calc\big[|\tilde{b}_t|^p\big]+E_\calc\big[|a_t|^p\big]+E_\calc\big[|\tilde{a}_t|^p\big]\big)
&\leq& K_p%\infty\quad a.s. 
\eeas
for every $p>1$.\noindent
\footnote{More precisely, the processes $\tilde{b}=(b_t)_{t\in[0,T]}$, $a=(a_t)_{t\in[0,T]}$, $\tilde{a}=(\tilde{a}_t)_{t\in[0,T]}$ and $\tilde{w}=(\tilde{w}_t)_{t\in[0,T]}$ are measurable mappings defined on $(\Omega,\calf)$, and for each $\omega'\in\Omega'$, 
the processes $\tilde{b}(\omega',\cdot)=(b_t(\omega',\cdot))_{t\in[0,T]}$, $a(\omega',\cdot)=(a_t(\omega',\cdot))_{t\in[0,T]}$, $\tilde{a}(\omega',\cdot)=(\tilde{a}_t(\omega',\cdot))_{t\in[0,T]}$ and $\tilde{w}(\omega',\cdot)=(\tilde{w}_t(\omega',\cdot))_{t\in[0,T]}$ 
satisfy the required conditions. 
$\tilde{w}(\omega',\cdot)$ is independent of $w(\omega',\cdot)$, 
i.e., $\tilde{w}$ and $w$ are $\calc$-conditionally independent, 
though this independency is not indispensable. 
}
\ed
\ed

Let 
\beas 
\bbY(\theta) &=& 
-\frac{1}{2T}\int_0^T \bigg\{\log\frac{\det S(\gamma,X_t,\theta)}{\det S(\gamma,X_t,\theta^*)}
+\text{Tr}\bigg(S(\gamma,X_t,\theta)^{-1}S(\gamma,X_t,\theta^*)-I_\sfm\bigg)\bigg\}dt.
\eeas
Define $\chi_0$ by 
\beas 
\chi_0 &=& \inf_{\theta\not=\theta^*}\frac{-\bbY(\theta)}{|\theta-\theta^*|^2}.
\eeas

\bd
\im[[D2\!\!]] For every $L>0$, 
\beas 
\sup_{r\in\bbN}r^LP_\calc\big[\chi_0\leq r^{-1}\big] &<& \infty\quad a.s.
\eeas
\ed

\begin{remark}\rm
If $\sigma$ does not depend on $\gamma$, then 
estimation with $\bbH_n$ needs no information about $\gamma$. 
\begin{en-text}
We assumed $\gamma$ was observable. 
If the function $\sigma$ depends only on $(x,\theta)$, then 
$\gamma$ can be unobservable, while $X_t$ may depend on $\gamma$ implicitly. 
In that case, the random field $\bbH_n$ needs no information about $\gamma$. 
\end{en-text}
\end{remark}

\begin{remark}\rm
We do not assume unconditional $L^p$ integrability of the functionals. 
For example, consider 
$X_t = \int_0^t e^{B_s^4}X_sds+\tilde{w}_t$ and 
the diffusion coefficient $\sigma(\gamma,X_t,\theta)=\theta \sqrt{1+X_t^2} $ for a Wiener process $B=(B_t)_{t\in[0,T]}$
living in $\calc$. 
Then $\sigma(\gamma,X_t,\theta)$ 
is not integrable. 
This situation is not formally treated in Uchida and Yoshida \cite{uchida2013quasi}. 
We can obtain a limit theorem for the QBE even in such a case. 
\end{remark}

\begin{remark}\rm
An analytic criterion and a geometric criterion for Condition $[D2]$ are provided by  \cite{uchida2013quasi}. %These sufficient conditions are easy to verify. 
\end{remark}

The random matrix $\Gamma$ is defined by 
\beas 
\Gamma &=& 
\frac{1}{2T}\int_0^T\text{Tr}\big((\partial_\theta S)S^{-1}(\partial_\theta S)S^{-1}(\gamma,X_t,\theta^*)\big)dt.
\eeas
We are writing 
$\hat{\theta}_n^M=\hat{\theta}_T$ and $\hat{\theta}_n^B=\tilde{\theta}_T$, and also 
$\hat{u}_n^M=\hat{u}_T$ and $\hat{u}_n^B=\tilde{u}_T$. 
We consider an extension $(\overline{\Omega},\overline{\calf},\overline{P})$ of $(\Omega,\calf,P)$. 
$\zeta$ denotes a random vector defined on this extension, 
having the $\sfp$-dimensional standard normal distribution $N_\sfp(0,I_\sfp)$ independent of $\calf$. 

\begin{theorem}
Suppose that $[D1]$ and $[D2]$ are fulfilled. Then, for every ${\sf A}\in\{M,B\}$, it holds that 
\beas 
E_\calc\big[f(\hat{u}_T^{\sf A})Y\big] &\to^P& \overline{E}_\calc\big[f(\Gamma^{-1/2}\zeta)Y\big]
\eeas
as $T\to\infty$ for any $\calf$-measurable bounded random variable $Y$ and any $f\in C(\bbR^\sfp)$ 
of at most polynomial growth. 
In particular, $\hat{u}_T^{\sf A}\to^{d_s(\calf)}\Gamma^{-1/2}\zeta$ as $T\to\infty$ for ${\sf A}=M,B$. 
\end{theorem}
\proof This result can be proved if we follow the proof of Theorems 4 and 5 of Uchida and Yoshida \cite{uchida2013quasi} 
in their Section 8, with the expectation $E$ replaced by the conditional expectation $E_\calc$. We omit details. 
$L_p$-boundedness of functionals are necessary, but it is possible under $E_\calc$ since the semimartingale structure 
is assumed under each $\calb_{\omega'}$.
\qed\halflineskip

\begin{remark}\label{290920-2}\rm 
Seemingly, we only considered time-independent scenario of the random field $\sigma$ 
represented by $\gamma$. % to simplify presentation. 
However, it is possible to consider a time-dependent coefficient $\sigma(t,\gamma,X_t,\theta)$ 
if we take $(t,X_t)$ for $X_t$. 
Then, this model includes also the model $\sigma(t,\gamma_t,X_t,\theta)$ 
having a time-varying component $\gamma=(\gamma_t)$. 
If we only assume discrete time observations $(\gamma_{t_j})$ of $\gamma$, then some condition for continuity of $\gamma$ would give  
similar results for the estimators. 
\end{remark}

%\bibliographystyle{plain}
%\bibliography{bibtex20080401}
\bibliographystyle{spmpsci}      % mathematics and physical sciences
\bibliography{bibtex-20160926-20170811-20170923}   % name your BibTeX data base

\end{document}
%%%%%%%%%%%%%%%%%%%%%%%%%%%%%%%%%%%%%%
%%%%%%%%%%%%%%%%%%%%%%%%%%%%%%%%%%%%%%
%%%%%%%%%%%%%%%%%%%%%%%%%%%%%%%%%%%%%%
%%%%%%%%%%%%%%%%%%%%%%%%%%%%%%%%%%%%%%
%%%%%%%%%%%%%%%%%%%%%%%%%%%%%%%%%%%%%%

\section{Introduction}

Yoshida \cite{Yoshida1997}, 
Uchida and Yoshida \cite{UchidaYoshida2015sVIC}

$\ep$ $\half$
${\bm A} {\bm \Phi}$
{\colorr a}{\coloroy b}{\colorr c}{\colorb d}
$\dotc$$\dot{C}$

\begin{theorem*}\label{th-1}
\bea\label{eq-1} 
a=b
\eea
\end{theorem*}
Thorem \ref{th-1} gives 
\begin{corollary*}
$b=c$
\end{corollary*}
(\ref{eq-1}): $a=b$

\section{Results}
\begin{theorem*}\label{th-2}
\bea\label{eq-1} 
a=b
\eea
\end{theorem*}
Thorem \ref{th-1} gives 
\begin{corollary*}
$b=c$
\end{corollary*}
(\ref{eq-1})

\begin{comment}
asdf
\end{comment}

%\bibliographystyle{plain}
%\bibliography{bibtex20080401}
\bibliographystyle{spmpsci}      % mathematics and physical sciences
\bibliography{bibtex-20150215}   % name your BibTeX data base

\end{document}